\documentclass[final,1p,times]{elsarticle}

\usepackage{amssymb}
\usepackage{amsmath}
\usepackage{graphicx}
\usepackage{xspace}
\usepackage{soul}  % ten pakiet pozwala uzywac komendy \st{} do przekreslania tekstu
\usepackage[ruled]{algorithm2e}
\usepackage{multirow}

\journal{Communications in Nonlinear Science and Numerical Simulation}

\newcommand{\capd}{CAPD::DynSys\xspace}
\newcommand{\PM}{\mathcal{P}}
\newcommand{\Real}{\mathbb{R}}

\usepackage{xcolor,listings}
\newtheorem{theorem}{Theorem}
\newtheorem{lemma}[theorem]{Lemma}
\newtheorem{remark}[theorem]{Remark}
\newtheorem{definition}[theorem]{Definition}
\newtheorem{corollary}[theorem]{Corollary}
\newenvironment{proof}{\textbf{Proof:}}{\qed}

\newcommand{\dom }{\mathrm{dom}\,}
\newcommand{\dt}{\quad:\quad}
\newcommand{\scprod}[2]{\langle #1\,|\, #2 \rangle}
\usepackage{color}

\input{intv_def.sty}

\begin{document}
    \lstset{language=C++,
        keywordstyle=\color{blue}\bfseries,
        commentstyle=\color{gray},
        stringstyle=\ttfamily\color{red!50!brown},
        showstringspaces=false}

    \begin{frontmatter}
        \title{Recent advances in rigorous computation of Poincar\'e maps}

        \author{Tomasz Kapela\fnref{MMMaestro}}
        \ead{Tomasz.Kapela@uj.edu.pl}

        \author{Daniel Wilczak\fnref{MMMaestro}\corref{ca}}
        \ead{Daniel.Wilczak@ii.uj.edu.pl}

        \author{Piotr Zgliczy\'nski\fnref{PZ}\corref{}}
        \ead{Piotr.Zgliczynski@ii.uj.edu.pl}

        \address{Faculty of Mathematics and Computer Science,
            Jagiellonian University, \L ojasiewicza 6, 30-348 Krak\'ow, Poland.}
        \fntext[MMMaestro] {This research is partially supported by
            the Polish National Science Center under Maestro Grant No. 2014/14/A/ST1/00453.}
        \fntext[PZ] {This research is partially supported by
            the Polish National Science Center under Grants No. 2019/35/B/ST1/00655 and UMO-2016/22/A/ST1/00077.}

        \cortext[ca] {Corresponding author.}
        \date{\today}
        \begin{abstract}
In this article we present recent advances on interval methods for rigorous computation of Poin\-car\'e maps.
We also discuss the impact of choice of Poincar\'e section and coordinate system on obtained bounds for computing Poincar\'e map nearby fixed points.
        \end{abstract}

        \begin{keyword}
            Poincar\'e map\sep
            rigorous numerical analysis\sep
            computer-assisted proof


            \sep % NA.Algorithms with automatic result verification
             % Periodic orbits of vector fields and flows
            %code \sep code (2000 is the default)
        \end{keyword}
    \end{frontmatter}

\section{Introduction}

Let us consider an ODE
\begin{equation}
 \label{eq:ode}
 x'(t)=f(x(t)),
\end{equation}
where $f\colon \mathbb R^n\to \mathbb R^n$ is a smooth (in most cases analytic and 'programmable')
function, and a local flow $\varphi(t,x)$ induced by (\ref{eq:ode}). The dynamics of continuous time systems can be efficiently studied by reduction to discrete time systems by introducing Poincar\'e sections and associated Poincar\'e maps. Such an approach is proved to be very efficient in studying periodic orbits, chaotic horseshoes, chaotic attractors, invariant tori and connecting orbits \cite{BarrioMartinezSerranoWilczak2015,CyrankaWanner,KokubuWilczakZgliczynski2007,SzczelinaZgliczynski2013,Wilczak2006,Wilczak2009abundance,WilczakZgliczynski2009siads,Tucker2002,TU2002127,Wilczak2010,WilczakBarrio2017,BarrioWilczak2020,Zelawski1999} or analysis of hybrid systems \cite{Galias2005, Galias2020,ISHII2009144,KONECNY20161}.

Informally, a Poincar\'e map $\mathcal P\colon \Pi_1\to \Pi_2$ assigns to a point $x\in \Pi_1$ the first intersection of its trajectory with $\Pi_2$, provided it exists --- precise definitions will be given in Section \ref{sec:notation}. There are various rigorous ODE solvers \cite{CAPD,CAPDREVIEW,Nedialkov2006,BerzMakino1999,KV,RauhBrillGunter2009} which compute enclosures of the solutions to IVPs. Based on them, one can implement an algorithm for computation of $\mathcal P(X)$, where $X\subset \Pi_1$. To the best of our knowledge, the \capd library \cite{CAPD,CAPDREVIEW} is the only freely available software that directly supports computation of enclosures of Poincar\'e maps.

The aim of this paper is to present recent advances in computation of Poincar\'e maps already implemented in \capd, for which the authors of this article are the main developers. We would like to discuss various aspects that impact obtained enclosures and reduce unavoidable overestimation. In particular, we would like to address the following questions.
\begin{itemize}
    \item What is the optimal choice (in the sense of obtained enclosures) of a Poincar\'e section near a periodic orbit?
    \item If the sections are fixed (like guards in hybrid systems), can we reduce overestimation by appropriate choice of coordinates on these sections?
    \item Can we take advantage from the knowledge on internal representation of the solution to IVP in a rigorous ODE solver to reduce overestimation in computation of Poincar\'e maps?
\end{itemize}

The article is organized as follows. In Section~\ref{sec:algorithm} we discuss our main algorithm and in Section~\ref{sec:choosingSection} we
discuss the questions listed above and report on multiple numerical experiments.

\subsection{Notation and basic definitions.}
\label{sec:notation}
An interval arithmetic is a base for rigorous numerical methods described in this paper.
We denote intervals and vectors and matrices with interval entries using bold face letters e.g. $\x$, $\X$. The interval hull of a set $S$, that is the smallest Cartesian product of closed intervals that contains $S$, is denoted by $\ihull{S}$. For a function $f\colon \mathbb R^n\to\mathbb R^m$ by $\fhull{f}$ we denote its (fixed) realisation in interval arithmetics, that is for $X\subset \mathbb R^n$ there holds $f(X)\subset \fhull{f}(X)$ provided the left hand side exists.

In what follows by $\varphi(t,x)$ we denote a local flow induced by (\ref{eq:ode}) with smooth vector field $f$.

For a matrix $A$ by $A^T$ we will denote the transposition of $A$.

For a submanifold $M$ of $\mathbb{R}^n$ and $y \in M$ by $T_y M \subset \mathbb{R}^n$ we will denote the tangent space to $M$ at $y$.

Throughout the paper for a scalar function $\alpha: \mathbb{R}^n \to \mathbb{R}$ the derivative $D\alpha (x)$ is the $1$-form on $T_x \mathbb{R}^n=\mathbb{R}^n$. It can be represented as
 \begin{equation*}
 D \alpha(x)=\left(\frac{\partial \alpha}{\partial x_1}(x),\dots,\frac{\partial \alpha}{\partial x_1}(x) \right).
 \end{equation*}
    This is a "row vector" as opposed of $\nabla \alpha(x)=(D \alpha (x))^T$, which is a "column vector". In fact the notion of $\nabla \alpha(x)$ depends on the scalar product,  which is silently assumed to be $\scprod{x}{y}=\sum_i x_iy_i$, hence is not independent from the coordinate system.

To denote the application of $1$-form $w$ to a vector $v$ we will use the following notation
\begin{equation}
  w(v)=\scprod{w}{v}.
\end{equation}

Usually in the text we will write vectors as $x=(x_1,\dots,x_n)$, i.e. as a row. However in the formulas we will  carefully distinguish between vectors as columns and $1$-forms as rows.

A $1$-form $w: \mathbb{R}^n \to \mathbb{R}$, $w\nequiv 0$   will be called a left-eigenvector of a matrix $A \in \mathbb{R}^{n \times n}$ if there exists $\lambda \in \mathbb{C}$, such that
\begin{equation}
  A^* w = \lambda w, \label{eq:leigenvect}
\end{equation}
where $A^*$ denotes the adjoint operator. When $w=(w_1,\dots,w_n)$ in the dual base, then (\ref{eq:leigenvect}) is equivalent to
\begin{equation}
  A^T w^T=\lambda w^T.
\end{equation}
It is well known that $A^*$ and $A$ have the same spectrum and if $\lambda$ is an eigenvalue of multiplicity $1$ with corresponding eigenvector $v$
and the left-eigenvector $w$, then
\begin{equation}
 \scprod{w}{v} \neq 0.
\end{equation}
Moreover, if $w_1$ is a left-eigenvector with eigenvalue $\lambda_1$ and $v_2$ is an eigenvector with eigenvalue $\lambda_2$, and $\lambda_1 \neq \lambda_2$, then
\begin{equation}
  \scprod{w_1}{v_2}=0.
\end{equation}

\begin{definition}
	The set $\Pi\subset \mathbb R^n$ is called a \emph{Poincar\'e section} or simply a \emph{section} for (\ref{eq:ode}) if $\Pi$ is connected manifold of codimension 1 and at every point $x\in\Pi$ vector field f is transversal to $\Pi$.
\end{definition}

In practice it is very convenient to define a Poincar\'e section $\Pi$ as a subset of zeroes of a smooth scalar function
\begin{equation}\label{eq:sectionFunction}
	\alpha\colon\mathbb R^n\to \mathbb R
\end{equation}
for which zero is a regular value. In the sequel we will always assume that a section is given by
\begin{equation}\label{eq:section}
	\Pi  = \Pi_{\alpha,\mathcal C} = \{x:\alpha(x)=0\wedge D\alpha(x)(f(x))\neq 0 \wedge \mathcal C(x)\},
\end{equation}
where $\mathcal C$ is a predicate that stands for additional constraints on the section (like crossing direction, restriction on the domain, etc.). The transversality condition $D\alpha(x)(f(x))\neq 0$ means that the flow is not tangent to $\Pi$ at every $x\in \Pi$. We will often write $\Pi$ or $\Pi_\alpha$ instead of $\Pi_{\alpha,\mathcal C}$ if $\alpha$ and/or $\mathcal C$ are clear from the context or irrelevant for the reasoning.

\begin{definition}
Let $\Pi$ be a section. We define a map $t_\Pi:\mathbb R^n\to \mathbb R$ in the following way
 \begin{enumerate}
  \item $x\in\dom t_\Pi$ iff there is $t>0$ such that $\varphi(t,x)\in \Pi$,
  \item for $x\in \dom t_\Pi$ we set $t_{\mathcal P}(x) = \inf\left\{t>0 : \varphi(t,x)\in \Pi\right\}$.
 \end{enumerate}
 We will call the function $t_\Pi$ the \emph{flow time} to section $\Pi$.
\end{definition}

\begin{definition}\label{def:PoincareMap}
 Let $\Pi_1,\Pi_2$ be sections for (\ref{eq:ode}). A Poincar\'e map $\mathcal P:\Pi_1\to\Pi_2$ is defined by $\mathcal P(x)= \varphi(t_{\Pi_2}(x),x)$ provided $t_{\Pi_2}(x)$ exists. We will denote the restriction of $t_{\Pi_2}$ to $\Pi_1$ by $t_{\mathcal P}$.
\end{definition}

\section{The main algorithm}
\label{sec:algorithm}
In this section we present the main algorithm for computation of rigorous enclosures of Poincar\'e maps. Here we assume that Poincar\'e sections $\Pi_1,\Pi_2$ and associated Poincar\'e map $\mathcal P:\Pi_1\to \Pi_2$ are fixed. We also assume that an affine change of coordinates $x\to A(x-y)$ is given, where $A$ is a linear isomorphism and $y\in \Pi_2$ is a vector.

Our goal is to construct an algorithm that for given $X\subset \Pi_1$ computes as tiny as possible enclosure $Y$ such that
\begin{equation}
    A(\mathcal P(X)-y)\subset Y.  \label{eq:ma-cc-sub}
\end{equation}
Let us comment briefly on the role of the choice of the above coordinate system, as it crucial in further analysis. The set $\mathcal{P}(X)$ is treated as a subset of $\mathbb{R}^n$. Although in our algorithm sections can be nonlinear, the algorithm was motivated by the simplest case, when $\Pi_2$ is a hyperplane. Assume $y\in \Pi_2$. The matrix $A$ should be chosen as $A=B^{-1}$, where one column (say first) of $B$ is transverse to $\Pi_2$
%and it defines direction of projection of $\mathbb R^n$ onto section $\Pi_2$,
while the remaining columns of $B$ span $T_y\Pi_2$, that is they give a coordinate system $(0,x_2,\ldots,x_n)\to y+B(0,x_2,\ldots,x_n)$ on $\Pi_2$. With this choice of $A$, computation of an enclosure of $\PM(X)$ in coordinates $(0,x_2,\ldots,x_n)$ is simply projection of the set $\varphi(t_{\Pi_2}(X),X)$ onto all coordinates but the first. The first column of $B$ controls direction of projection of $\mathbb R^n$ onto section $\Pi_2$, while proper choice of the basis of $T_y\Pi_2$ helps in reducing wrapping effect on the section.

The algorithm consists of the following two steps. First, given a rough enclosure for the return time $t_{\Pi_2}(X)$ we try to obtain sharp bounds for it --- Section~\ref{sec:returnTime}. Then, we compute $\PM(X)$ in a given affine coordinate system, that is an enclosure $Y$ for $A(\mathcal P(X)-y)$  --- Section~\ref{sec:computePoincareMap}.
To suppress overestimation in both steps we take advantage from the knowledge about internal representation of subsets of $\mathbb R^n$ used by the underlying rigorous ODE solver. Abstract assumptions on these data structures and their sample realisations will be presented in Section \ref{sec:representation}.

\subsection{Representation of sets in rigorous ODE solvers}\label{sec:representation}
Every rigorous ODE solver uses its own internal representation of subsets of $\mathbb R^n$ and solutions to IVP. When a set of initial condition is propagated by a dynamical system and on each step the image is bounded by an interval vector (product of intervals), then typically we observe the wrapping effect that leads to huge overestimation. On the other hand when the image is bounded by some non-linear shape, e.g given by multidimensional polynomials (like in the case of Taylor models \cite{BerzMakino1999}), then the result is more accurate but the computational cost increases rapidly with the dimension and the degree of the polynomial.

In \cite{CAPD,CAPDREVIEW,Nedialkov2006} the sets are represented (see \cite{MrozekZgliczynski2000}) as parallelepipeds, doubletons and tripletons. These strategies provide a good compromise between speed and accuracy as shown in \cite{MIYAJI201634}.

A \emph{doubleton} is a tuple $(x,C,\rz,Q,\q)$, where $x\in\mathbb R^n$, $\rz,\q\subset\mathbb R^n$ are interval vectors such that $0\in \rz\cap \q$ and $C,Q\in\mathbb R^{n\times n}$ are matrices with $Q$ being invertible and usually close to orthogonal. It represents set of points in $\mathbb R^n$ defined by
\begin{equation}\label{eq:doubleton}
X = \left\lbrace x + C r_0 + Q q \mid r_0 \in \rz, q \in \q \right\rbrace.
\end{equation}

A \emph{tripleton} is defined as a tuple $(x,C,\rz,Q,\q,B,\r)$  where $x\in\mathbb R^n$, $\rz,\q,\r\subset\mathbb R^n$ are interval vectors containing 0 and $C,Q,B\in\mathbb R^{n\times n}$ are matrices. Geometrically it is an intersection of two doubletons, that is
\begin{equation}\label{eq:tripleton}
X = \left\lbrace x + C r_0 + Q q \mid r_0 \in \rz, q \in \q \right\rbrace \cap \left\lbrace x + C r_0 + B r \mid r_0 \in \rz, r \in \r \right\rbrace.
\end{equation}
The matrix $Q$ is assumed to be close to orthogonal and invertible, while on $B$ we only assume it is invertible.

In the sequel we assume that an abstract data structure called \texttt{RepresentableSet} is defined and it represents a subset of $\mathbb R^n$. For example as \texttt{RepresebtableSet} one can take doubleton or tripleton. We also assume that for given \texttt{RepresentableSet} two algorithms \texttt{eval} and \texttt{affineTransform} are provided. Their expected  inputs and outputs are defined by Algorithm~\ref{alg:eval} and Algorithm~\ref{alg:affineTransform}. For doubleton representation their possible realisations are given by Algorithm~\ref{alg:doubletonEval} and Algorithm~\ref{alg:doubletonAffineTransform}, respectively.

\begin{algorithm}\label{alg:eval}
    \KwIn{$X\subset \mathbb R^n$ \dt \texttt{RepresentableSet}}
    \KwIn{$g:\mathbb R^n\to \mathbb R^m$ \dt smooth function}
    \KwOut{An enclosure for $g(X)$}
    \caption{{\sc eval}}
\end{algorithm}
\begin{algorithm}\label{alg:affineTransform}
    \KwIn{$X\subset \mathbb R^n$ \dt \texttt{RepresentableSet}}
    \KwIn{$A:\mathbb R^n\to \mathbb R^m$ \dt linear map}
    \KwIn{$y\in \mathbb R^n$}
    \KwOut{An enclosure of $A(X-y)$}
    \caption{{\sc affineTransform}}
\end{algorithm}
\begin{algorithm}
    \KwIn{$x + C\rz + Q\q\subset\mathbb R^n$ \dt \texttt{Doubleton}}
    \KwIn{$g:\mathbb R^n\to \mathbb R^m$ \dt smooth function}
    \KwOut{An enclosure of $g(x+C\rz+Q\q)$}
    \BlankLine
    $\X\gets \ihull{x + C\rz + Q\q}$; \qquad // enclose set as interval vector \\
    $\M\gets \fhull{Dg}(\X)$; \qquad // enclose derivative as interval matrix \\
    \Return{$\fhull{g}(\X) \cap \ihull{\fhull{g}(x) + (\M C)\rz + (\M Q)\q}$};
    \caption{{\sc eval}}
    \label{alg:doubletonEval}
\end{algorithm}
\begin{algorithm}[htbp]
    \KwIn{$x + C\rz + Q\q \subset\mathbb R^n$  \dt \texttt{Doubleton}}
    \KwIn{$A:\mathbb R^n\to \mathbb R^m$ \dt linear map}
    \KwIn{$y\in \mathbb R^n$}
    \KwOut{An enclosure for $A( x+C\rz + Q\q-y)$}
    \BlankLine
    \Return{$\ihull{A((x-y) + C\rz + B\r)} \cap \ihull{A(x-y) + (AC)\rz + (AQ)\q}$}\;
    \caption{{\sc affineTransform}}
    \label{alg:doubletonAffineTransform}
\end{algorithm}

Observe, that the assumption $0\in \rz\cap \q$ guarantees that $x\in X$ and thus \\
$\ihull{\fhull{g}(x) + (MC)\rz + (MQ)\q}$ is an enclosure of $g(x+C\rz+Q\q)$ in Algorithm~\ref{alg:doubletonEval}.

\subsection{Refinement of return time $t_{\mathcal P}(X)$}\label{sec:returnTime}
Assume that sections $\Pi_1, \Pi_2 \subset \Real^n$ are defined as sets of zeroes of  smooth functions $\alpha_1, \alpha_2: \Real^n \to \Real$ .
Computation of a rough enclosure of the return time $t_{\Pi_2}(X)$ for $X \in \Pi_1$ can be realized in the following steps:
\begin{enumerate}
 \item Find $t_1>0$ such that $t_{\mathcal P}(X)\geq t_1$. This can be done by checking $\varphi([\delta,t_1],X)\cap \Pi_2=\emptyset$ where $\delta$ is the lower bound for the return time $t_{\mathcal P}(X)$ coming from transversality of the flow for both sections.
 \item Find $t_2>t_1$ such that for $x\in X$ there holds $\alpha_2(\varphi(t_2,x))\alpha_2(\varphi(t_1,x))< 0$. This guarantees that the function $\alpha_2$ changes the sign on the solution segment $\varphi([t_1,t_2],x)$.
 \item Check transversality condition:
 \begin{equation}\label{eq:transversalityCondition}
 D\alpha_2(\varphi(t,x))\left(f(\varphi(t,x))\right)\neq 0 \text{ for } t \in [t_1,t_2], x \in X.
 \end{equation}
\end{enumerate}
If all steps do not fail then $[t_1,t_2]$ is an upper bound for the return time $t_{\mathcal P}(X)$. This rough estimation can be significantly improved by applying Interval Newton Operator \cite{Moore1966} to the following univariate function
\begin{equation}\label{eq:returnTime}
 g_x(t):=\alpha_2(\varphi(t,x)).
\end{equation}
\begin{lemma}[Interval Newton Operator]
Let $g:\mathbb R^n\to\mathbb R^n$ be a smooth function and let $X$ be a convex, closed set $X$. For fixed $x_0\in X$ the interval Newton operator is defined as
\begin{equation}\label{eq:INO}
N(g,x_0,X) = x_0 - \ihull{Dg(X)}^{-1}\fhull{g}(x_0).
\end{equation}
\begin{itemize}
    \item If $N(g,x_0,X)\subset \mathrm{int}(X)$ then the $g$ has in $X$ unique zero $x_*$ such that $x_*\in N(g,x_0,X)$.
    \item If $x\in X$ is a zero of $g$ then $x\in N(g,x_0,X)$.
\end{itemize}
\end{lemma}

\begin{algorithm}
\KwIn{$[t_1,t_2]$ \dt an interval that encloses $t_{\Pi_2}(X)$}
\KwIn{$X_1$ \dt \texttt{RepresentableSet} that encloses $\varphi(t_1,X)$}
\KwIn{$\alpha_2$ \dt function that defines section $\Pi_2$}
\KwIn{$f$ \dt underlying vector field}
\KwOut{Improved bound for $t_{\Pi_2}(X)$}
\BlankLine
$t_0\gets (t_1+t_2)/2$\;
$X_0\gets$ \texttt{RepresentableSet} that encloses $\varphi(t_0-t_1,X_1)$;\qquad // from ODE solver\\
$\intv{g_0}\gets \texttt{eval}(X_0,\alpha_2)$\;
$\intv{e}\gets \texttt{eval}(X_1,\varphi([0,t_2-t_1],\cdot))$;\quad // enclosure from ODE solver\\
$\intv{g}\gets \texttt{eval}\left(\intv{e},\nabla\alpha_2(\cdot)\cdot f(\cdot)\right)$\;
\Return{$[t_1,t_2]\cap (t_0-\intv{g_0}/\intv{g})$}\;
\caption{{\sc refineReturnTime}}
\label{alg:refineReturnTime}
\end{algorithm}

\begin{lemma}
Assume the Algorithm~\ref{alg:refineReturnTime} is called with the correct arguments and returns $T$. Then $ t_{\Pi_2}(X)\subset T\subset[t_1,t_2]$.
\end{lemma}
\begin{proof}
Let us fix $x\in X$. By the assumptions of the Algorithm~\ref{alg:refineReturnTime} we have $x_1:=\varphi(t_1,x)\in X_1$.
Moreover, since $t_{\Pi_2}(x)\in[t_1,t_2]$ we have $t_{\Pi_2}(x) = t_1 + t_{\Pi_2}(x_1)$. We will apply the Interval Newton Operator \cite{Moore1966} to the function $g_{x_1}$ defined by (\ref{eq:returnTime}) in order to obtain a tighter enclosure for
$$
t_{\Pi_2}(x_1)\in [0,t_2-t_1]\cap g_{x_1}^{-1}(0).
$$
The interval Newton operator applied to $g_{x_1}$ reads
\begin{equation*}
 N(g_{x_1},t_0-t_1,[0,t_2-t_1]) = (t_0-t_1) - \fhull{g_{x_1}}(t_0-t_1)/\fhull{g_{x_1}'}([0,t_2-t_1]).
\end{equation*}
It is easy to see that for $x\in X$ the algorithm encloses the following quantities with $t_0=(t_1+t_2)/2$ and $t\in[0,t_2-t_1]$
\begin{eqnarray*}
 g_{x_1}(t_0-t_1)&=& \alpha_2(\varphi(t_0-t_1,x_1))\in \intv{g_0},\\
 g_{x_1}'(t) &=& \nabla\alpha_2(\varphi(t,x_1))\cdot f(\varphi(t,x_1)) \in \intv{g}.
\end{eqnarray*}
Therefore $t_{\Pi_2}(x_1)\in (t_0-t_1)-\intv{g_0}/\intv{g}$ and in consequence
\begin{equation*}
 t_{\Pi_2}(x) = t_1 + t_{\Pi_2}(x_1) \in t_0 - \intv{g_0}/\intv{g},
\end{equation*}
which completes the proof.
\end{proof}

Observe that $\intv{g}$ does not contain zero because the transversality condition (\ref{eq:transversalityCondition}) is satisfied on $[t_1,t_2]$. Algorithm~\ref{alg:refineReturnTime} can be iterated until no improvement is observed. In most cases it stabilizes after $2-3$ iterations provided $[t_1,t_2]$ was already a good initial bound for $t_{\mathcal P}(X)$. Observe that $diam(t_{\mathcal P}(X))$ is proportional
to $diam(\alpha_2(X_0))$, which is up to a constant the size of $X_1$ in the "direction" of the section.

\begin{remark}
It is very important to use information about representation of the set $X_0$ in computation $\intv{g_0}\gets {\emph{\texttt{eval}}}(X_0,\alpha_2)$.
For example, for $\alpha_2(x,y)=x-y$ and the set
$X_0=x+C r=\begin{bmatrix}1+\varepsilon\\1\end{bmatrix}+\begin{bmatrix}1 & 1 \\ 1 & -1\end{bmatrix}\begin{bmatrix}[-1,1]\\0\end{bmatrix}$
a direct evaluation in interval arithmetics gives
$$
\alpha_2(X_0) \subset \alpha_2\left((1+\varepsilon,1) + ([-1,1],[-1,1])\right) = \alpha_2\left(([\varepsilon,2+\varepsilon],[0,2]\right) = [-2+\varepsilon,2+\varepsilon].
$$
The diameter of obtained estimate is 4 and as it contains 0 we cannot conclude if the set $X_0$ is before or after the section.
Representing $X_0$ as doubleton $(x, C, \rz, Q=Id, \q=0)$ and applying Algorithm \ref{alg:doubletonEval} gives
%multiplication of ``thin'' objects first gives
$$
\alpha_2(X_0)\subset [\alpha_2](1+\varepsilon,1) + \left(\nabla\alpha_2(X_0)^T\begin{bmatrix}1 & 1 \\ 1 & -1\end{bmatrix}\right)\begin{bmatrix}[-1,1]\\0\end{bmatrix} =
[\varepsilon,\varepsilon] + \begin{bmatrix} 0 & 2\end{bmatrix} \begin{bmatrix}[-1,1]\\0\end{bmatrix} = [\varepsilon,\varepsilon].
$$
This sharp bound, even for $\varepsilon$ close to zero, allows to verify relative location of the set $X_0$ and the section.
\end{remark}

\subsection{Computation of Poincar\'e map}\label{sec:computePoincareMap}
In this section we present Algorithm~\ref{alg:computePoincareMap}, which realizes the last step of computation of an enclosure of $A(\mathcal P(x)-y)$ for $x\in X$, where $y\in\Pi_2$ and $A$ is a linear mapping. That is, we represent $\mathcal P(x)$ in a local affine coordinate system near $y$. Let us emphasize, that the restriction to affine coordinates is reasonable keeping in mind that the underlying ODE solver represents subsets of $\mathbb R^n$ usually as doubletons or tripletons.

\begin{algorithm}%[htbp]
\KwIn{$[t_1,t_2]$ \dt an interval that encloses $t_{\mathcal P}(X)$}
\KwIn{$X_1$ \dt \texttt{RepresentableSet} that encloses $\varphi(t_1,X)$}
\KwIn{$f$ \dt underlying vector field}
\KwIn{$y$ \dt a vector}
\KwIn{$A$ \dt a linear map}
\KwOut{An enclosure of $A(\mathcal P(X)-y)$}
\BlankLine

$\intv{e}\gets \texttt{eval}(X_1,\varphi([0,t_2-t_1],\cdot))$;\quad // enclosure from ODE solver\\
$t_0\gets (t_1+t_2)/2$\;
$\intv{\Delta t}\gets[t_1,t_2]-t_0$\;
$X_0\gets$ \texttt{RepresentableSet} that encloses $\varphi(t_0-t_1,X_1)$;\quad // call to ODE solver\\
$\intv{y_0}\gets \texttt{affineTransform}(X_0,A,y)$\;
$\intv{y}\gets \texttt{eval}(X_0,A\circ f)\cdot \intv{\Delta t}$\;
$\intv{\Delta y} \gets \frac{1}{2} A\cdot \fhull{Df}(\intv{e})\cdot \fhull{f}(\intv{e})\cdot \intv{\Delta t}^2$\;
$\intv{z}\gets (\intv{y_0}+\intv{y}+\intv{\Delta y})\cap \ihull{A(\intv{e}-y)}$\;
\Return{$\intv{z}$}\;
\caption{{\sc computePoincareMap}}
\label{alg:computePoincareMap}
\end{algorithm}	

Before we state and prove properties of Algorithm~\ref{alg:computePoincareMap}, let us comment again (see the paragraph after (\ref{eq:ma-cc-sub})) on the role of the coordinate system
\begin{equation*}
x \mapsto z=A(x-y),
\end{equation*}
as it crucial in further analysis. First, let us observe that the set $\mathcal P(X)\subset \Pi_2\subset \mathbb R^n$ is treated as a subset of full phase space. Formally the algorithm computes an enclosure for $\intv{z}=A(\mathcal P(X)-y)\subset \mathbb R^n$ but what we usually are interested in is indeed $\mathcal P(X)$. We will show, that the proper choice of the matrix $A$ helps to extract from computed result $A(\mathcal P(X)-y)$ much tighter enclosure on $\mathcal P(X)$ than from a naive algorithm, that is  $\mathcal P(X) \subset\texttt{eval}(X_1,\varphi([0,t_2-t_1],\cdot))$.

If the Poincar\'e section $\Pi_2$ is a hyperplane with $y\in \Pi_2$ then the matrix $A$ should be chosen as $A=B^{-1}$, where one column (say first) of $B$ is transverse to $\Pi_2$ and it defines direction of projection of $\mathbb R^n$ onto $\Pi_2$, while the remaining columns of $B$ should span $T_y\Pi_2$. Denote by $\intv{z}\supset A(\mathcal{P}(X)-y)$ the result of Algorithm~\ref{alg:computePoincareMap}. With this choice of $A$ we have
$$
\mathcal P(X) \subset y + B(0,\intv{z}_2,\ldots,\intv{z}_n)^T
$$
because all but first columns of $B$ span the section $T_y\Pi_2$. What we get is an enclosure of $\mathcal P(X)$ in the form of an affine set.

\begin{lemma}
\label{lem:algComputePoincareMap}
 Assume the Algorithm~\ref{alg:computePoincareMap} is called with its arguments, where $f$ is of class $\mathcal C^1$. Let $\intv{z}$ be computed result. Then $A(\mathcal P(X)-y)\subset \intv{z}$.
\end{lemma}
\begin{proof}
 First observe that the computed quantity $\intv{e}$ is an enclosure of $\mathcal P(X)$ and therefore $A(\mathcal P(X)-y)\subset A(\intv{e}-y)$. We have to show that $A(\mathcal P(X)-y)\subset \intv{y_0}+\intv{y}+\intv{\Delta y}$.

 Denote by $T=[0,t_2-t_1]$. The flow $\varphi$ is $\mathcal C^2$ with respect to time variable. Using the Taylor theorem with the Lagrange remainder we obtain
 \begin{equation*}
  \mathcal P(X)\subset \varphi(T,X_1)=\varphi(\intv{\Delta t},X_0)\subset X_0 + f(X_0)\intv{\Delta t} + \frac{1}{2}Df(\intv{e})f(\intv{e})\intv{\Delta t}^2,
 \end{equation*}
 where $X_0,\intv{e}, \intv{\Delta t}$ are defined as in Algorithm~\ref{alg:computePoincareMap}.
Therefore
\begin{eqnarray*}
A(\mathcal P(X)-y)&\subset& A(X_0-y) + Af(X_0)\intv{\Delta t} + \frac{1}{2}ADf(\intv{e})f(\intv{e})\intv{\Delta t}^2\\
&\subset& \intv{y_0}+\intv{y}+\intv{\Delta y}.
\end{eqnarray*}
\end{proof}
\begin{remark}
It is extremely important to compute $\intv{y_0}$ and $\intv{y}$ taking into account the representation of $X_0$ (Algorithms \emph{\texttt{eval}} and \emph{\texttt{affineTransform})} in order to reduce wrapping effect in evaluation of products $A\cdot X_0$ and $A\cdot f(X_0)$. These terms are the most important components of the result. The quantity $\intv{\Delta y}$ is expected to be small in comparison to $\intv{y_0}$. Indeed, if the diameter of the set is small, then we expect that the diameter of crossing time $[t_1,t_2]$, and so $\intv{\Delta t}$ is proportional to diameter of the set $X_0$ and thus $\intv{\Delta y} \in O\left(diam(X_0)^2\right)$.  More precisely $\Delta t$ is proportional to diameter of the set $X_0$ in the direction of the flow.
\end{remark}

In the next section we will argue, that playing with the choice of Poincar\'e section and coordinate system $A$ on it we can make the term $\intv{y}$ also small in comparison to $diam(\intv{y_0})$.

\section{Varying Poincar\'e section and coordinates on it.\label{sec:choosingSection}}

\begin{definition}
  Assume $p$ is a periodic point of (\ref{eq:ode}) of period $T$. $\lambda \in \mathbb{C}$ is \emph{a Floquet multiplier} for the periodic orbit through $p$ if it is an eigenvalue of $D_x\varphi(T,p)$.
\end{definition}
It is well known that the set of Floquet multipliers depends on the orbit, only, i.e. is the same for all choices of $p$ on the periodic orbit. Moreover, the set of Floquet multipliers always contains $\lambda=1$ corresponding to the orbit direction $f(p)$ and the eigenvalues for the Poincar\'e map on any transversal to $f(p)$ section containing $p$ give all Floquet multipliers different from $1$.

In the remainder of the section we will often say that a ``vector $f(p)$ is orthogonal to the section $\Pi$''. In saying that,  we silently assumed that we have a scalar product given by $\scprod{x}{y}=\sum_i x_iy_i$. Stated as such the notion of ``orthogonality'' depends on the coordinate frame used. The ``ortogonality'' in the above sense is used as a measure to assure that the vector field is far from being tangent to the section $\Pi$.

\begin{figure}[htbp]
\centerline{
    \includegraphics[width=.75\textwidth]{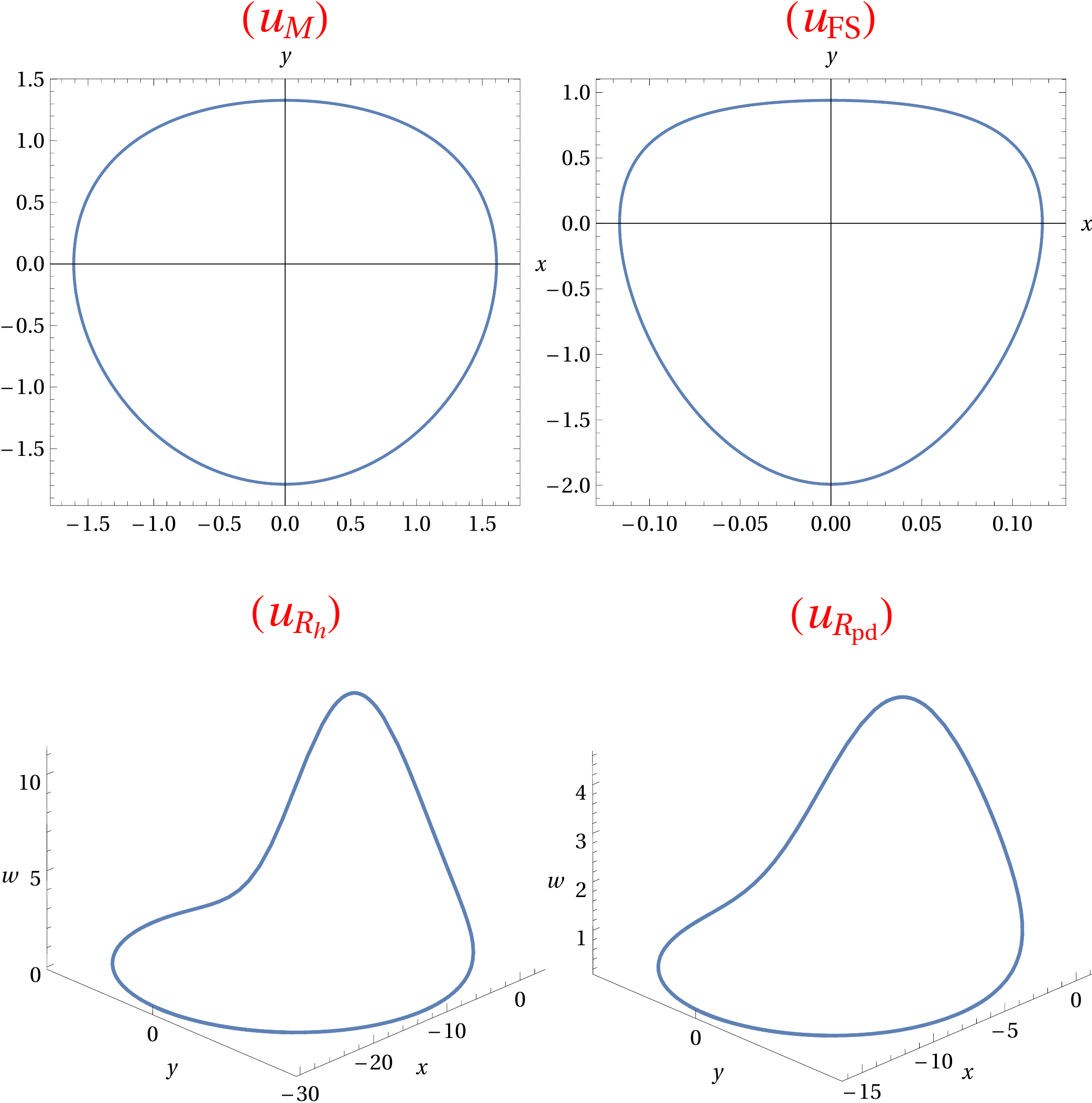}
}
\caption{The four periodic orbits $u_M$, $u_{FS}$, $u_{R_h}$ and $u_{R_{pd}}$, for the Michelson system \cite{Mi}, Falkner-Skan system \cite{FS} and the R\"ossler system \cite{Rossler79}, respectively.\label{fig:periodic-orbits}}
\end{figure}

As mentioned in Section~\ref{sec:computePoincareMap}, the term $\intv{y_0}$ is the main component of the result returned by Algorithm~\ref{alg:computePoincareMap} and computation of $\intv{y_0}$ takes into account representation of the sets by an underlying ODE solver (algorithm \texttt{affineTransform}). In this section we will give some heuristics on how one can reduce diameter of the component $\intv{y}$. Our considerations will be supported by the following case studies
\begin{itemize}
    \item the Michelson system \cite{Mi} --- a 3D system given by
    \begin{equation}\label{eq:michelson}x'=y,\quad y'=z,\quad z'=c^2-y-\frac{1}{2}x^2\end{equation}
    around an approximate periodic orbit
    \begin{equation}\label{eq:uM}u_{M}\approx(0,1.32825866108569290258,0)\end{equation}
    for parameter $c=0.8.$
    It has been shown \cite{Wilczak2006} that for a wide range of parameter values including $c=0.8$ certain Poincar\'e map of this system is $\Sigma_4$ chaotic and hyperbolic periodic orbit $u_M$ is part of chaotic invariant set. Approximate Floquet multipliers of $u_M$ are $1$ and
    \begin{equation}\label{eq:michelsonEV}\lambda_{M_1}\approx -21.57189303583905,\quad \lambda_{M_2}\approx -0.046356617768258279\end{equation}
    \item the Falkner-Skan equation \cite{FS} --- a 3D system given by
    \begin{equation}\label{eq:falknerskan}x'=y,\quad y'=z,\quad z'=c(y^2-1)-xz\end{equation}
    around an approximate periodic orbit
    \begin{equation}\label{eq:uFS}u_{FS}\approx(0,0.939712208779672476275,0)\end{equation}
    for parameter $c=250.$ As shown in \cite{WalawskaWilczak2019}, this periodic orbit belongs to a family parameterized by $c\in[\frac{9}{8},100000]$. The choice of $c=250$ is motivated by the fact, that this orbit appears to be ill--conditioned for the method of minimization of the crossing time $[t_1,t_2]$ --- see Section~\ref{sec:variableSection}. It is hyperbolic with approximate Floquet multipliers $1$ and
    \begin{equation}\label{eq:fsEV}\lambda_{FS_1}\approx -3.1255162015308575,\quad \lambda_{FS_2}\approx -0.31994714969329141.\end{equation}
    Moreover, the system is stiff near the initial condition, which may cause additional problems for an ODE solver.
    \item the R\"ossler system \cite{Rossler79} -- a 4D system given by
    \begin{equation}\label{eq:rossler4D}x'=-y-w,\quad y'=x+ay+z,\quad z'=dy+cw,\quad w'=xw+b\end{equation} around an approximate periodic orbit
    \begin{equation}\label{eq:uRH}u_{R_h}\approx (-29.841563300389689,0,15.047757539453583,0.10059818458161384)\end{equation} for parameter values $a=0.25$, $b=3$, $c=-0.5$, $d=0.05$. As shown in \cite{Wilczak2009abundance}, this system is hyperchaotic (two positive Lyapunov exponents) and $u_{R_h}$ belongs to hyperchoatic invariant set. Approximate Floquet multipliers are $1$ and
    \begin{equation}\label{eq:rosslerH}\lambda_{R_{h1}}\approx -2.9753618617897111,\quad \lambda_{R_{h2}}\approx 	1.11933293616997,\quad \lambda_{R_{h3}}\approx -2\cdot 10^{-18}\end{equation}
    \item the R\"ossler system (\ref{eq:rossler4D}) around an approximate periodic orbit
    \begin{equation}\label{eq:uRPD}u_{R_{pd}}\approx (-16.051468914417546,0,8.362179513564907,0.18738588995067224)\end{equation}
    for parameter values $a=0.25$, $b=3$, $c=-0.397617541005413$, $d=0.05$. This orbit belongs to the same family as $u_{R_h}$ parameterized by $c$. Near this particular parameter value the family undergoes the period doubling bifurcation. The Floquet multipliers are $1$ and
\begin{equation}\label{eq:rosslerPD}\lambda_{R_{pd1}}\approx 1.2039286263296654,\quad \lambda_{R_{pd2}}\approx -1,\quad \lambda_{R_{pd3}}\approx -6\cdot 10^{-17}\end{equation}
and absolute values of $\lambda_{R_{pd1}}$ and $\lambda_{R_{pd2}}$ are close or very close to one, hence non--linear terms can play significant role near this orbit. We expect that this orbit may be difficult for both ODE solver and our algorithm for computation of Poincar\'e map, because the ODE solver uses affine representation of sets and also the Algorithm~\ref{alg:computePoincareMap} computes an enclosure of $\mathcal P$ in affine form.
\end{itemize}
Shapes of these four periodic orbits are shown in Figure~\ref{fig:periodic-orbits}.

\subsection{Reduction of sliding and wrapping effects}\label{sec:fixedSections}
In this paragraph we assume that we do not have freedom in choosing the Poincar\'e section. This is common in hybrid systems (sections are guards between regions with different vector fields) or when one wants to take advantage of some specific properties of a system, like symmetries.

Recall, that $\intv{y_0}$ is the main component of the result returned by Algorithm~\ref{alg:computePoincareMap}. In what follows we will argue that the diameter of the projection of $\intv{y}$ onto section can be reduced by a proper choice of direction of projection, which is hidden in the coordinate system $A$.

\begin{theorem}\label{thm:sliding}
Assume the Algorithm~\ref{alg:computePoincareMap} is called with its arguments $y\in X_0\cap \Pi$ and the matrix $A=B^{-1}$, where
$$B = \left[\begin{array}{c|c}f(y) & M\end{array}\right]$$
and columns in $M$ span $T_y \Pi$. Then
$$\intv{y}+\intv{\Delta y}\in (\intv{\Delta t},0,0,\ldots) + O\left(diam(X_0)^2\right).$$
\end{theorem}
\begin{proof}
 In Algorithm~\ref{alg:computePoincareMap}, the quantity $\intv{y}$ is computed as $\intv{y}\gets \texttt{eval}(X_0,A\circ f)\cdot \intv{\Delta t}$ which is realized by the intersection of direct evaluation in interval arithmetic and application of the mean value theorem to $A\circ f$. Thus
    \begin{equation*}
    \intv{y} \subset \ihull{A f(x_0) + A \fhull{Df}(X_0) \Delta X_0}\intv{\Delta t},
    \end{equation*}
    where $x_0=\mathrm{mid}(X_0)$ and $\Delta X_0 = X_0-x_0$. We have
    \begin{equation}\label{eq:diamx01}
    A \fhull{Df}(X_0) \Delta X_0\intv{\Delta t}\in O\left(diam(X_0)^2\right)
    \end{equation}
    because $diam(\intv{\Delta t})\in O\left(diam(X_0)\right)$.
    Since $y\in X_0$, which is a convex set, we have
    $$
    A f(x_0)\intv{\Delta t} \in (Af(y))\intv{\Delta t} + \ihull{ADf(X_0)}(x_0-y)\intv{\Delta t}.
    $$
    By the assumptions on $A$, we have $Af(y)=(1,0,0,\ldots)$ and thus
    \begin{equation}\label{eq:diamx02}
    \ihull{(Af(y))\intv{\Delta} t} = (\intv{\Delta t},0,0\ldots).
    \end{equation}
    We also have
    \begin{equation}\label{eq:diamx03}
    \ihull{A Df(X_0)}(x_0-y)\intv{\Delta t}\in O\left(diam(\Delta X_0)diam(\intv{\Delta t})\right) = O\left(diam(\Delta X_0)^2\right).
    \end{equation}

Gathering (\ref{eq:diamx01})--(\ref{eq:diamx03}) we get $\intv{y} = (\intv{\Delta t},0,0\ldots) + O\left(diam(X_0)^2\right)$.
The component $\intv{\Delta y}$ contains the factor $\intv{\Delta t}^2$ and thus its diameter is also quadratic in diameter of $X_0$.
\end{proof}

\begin{remark}
    In the proof of Theorem~\ref{thm:sliding} we assumed that evaluations $f(x_0)$ and $f(y)$ are exact and thus of diameter zero.
\end{remark}

\begin{corollary} Assume that the section $\Pi$ is a hyperplane. Denote by $\intv{z}=\intv{y}+\intv{y_0}+\intv{\Delta y} = (\intv{z_1},\ldots,\intv{z_n})$ result of Algorithm~\ref{alg:computePoincareMap}. Under assumptions of Theorem~\ref{thm:sliding} and 
    $$
    \mathcal P(X) \subset y + B (0,\intv{z_2},\ldots,\intv{z_n})^T.
    $$
    Moreover,
    $$
        \mathcal P(X) \subset y + B (0,\intv{y_2},\ldots,\intv{y_n})^T + O\left(diam(X_0)^2\right).
    $$
\end{corollary}
\begin{proof}
    Since $y\in \Pi$ and all but first columns of $B$ span $T_y\Pi$, we have that
    $$\Pi =\left\{ y + B (0,z_2,\ldots,z_n)^T : z_2,\ldots,z_n\in\mathbb R\right\}.$$
    Therefore,
    $$\mathcal P(X) \subset (y+B\intv{z}) \cap \Pi = y + B (0,\intv{z_2},\ldots,\intv{z_n})^T.$$
    From Theorem~\ref{thm:sliding} we have that
    $$
        (0,\intv{z_2},\ldots,\intv{z_n}) = (0,\intv{y_2},\ldots,\intv{y_n}) + O\left(diam(X_0)^2\right),
    $$
    which completes the proof of the second assertion.
\end{proof}
\subsection{Choosing coordinate system on Poincar\'e section close to fixed point}
Let $u\in \Pi$ be an approximation of a fixed point of Poincar\'e Map $\PM: \Pi \to \Pi$ and assume $\Pi$ is an affine section
$$
\Pi=\left\{u+\sum_{i=2}^n z_iv_i\,|\,z_i\in\mathbb R \right\},
$$
where $\{v_2,\ldots,v_n\}$ are linearly independent. We introduce here our two strategies for setting up coordinate system centred at $u$. In the first step we compute non-rigorously $D\PM(u)$ in the coordinate system 
$$
(z_2,\ldots,z_n) = u + \sum_{i=2}^n z_iv_i.
$$
Let $M\in\mathbb{R}^{(n-1)\times(n-1)}$ be a square matrix with columns corresponding to normalized eigenvectors of $D\PM(u)$. In the second step we choose any vector $v_1$ linearly independent with $\{v_2,\ldots, v_n\}$ and define matrices  
\begin{equation}\label{eq:matricesVQB}
V = \left[v_1\ v_2\ \ldots v_n\right],\qquad 
Q = \begin{bmatrix}
	1 & 0\\
	0 & M
\end{bmatrix} \qquad\text{and}\qquad B = VQ.
\end{equation}
Depending on the choice of $v_1$ we call this strategy:
\begin{itemize}
  \item \textbf{diag+normal} if $v_1$ is orthogonal to the section $\Pi$,
  \item \textbf{diag+flowdir} if $v_1$ is normalized vector field at $u$ i.e. $v_1=\frac{f(u)}{\|f(u)\|}$.
\end{itemize}

\subsection{Choosing coordinate systems on Poincar\'e section --- experiments.}
In this section we present results of numerical experiments regarding various strategies of choosing the matrix $A$ in Algorithm~\ref{alg:computePoincareMap}. All results strongly confirm theoretical findings of Theorem~\ref{thm:sliding}. In our test suite, we have chosen ``standard'' sections, that is $x=0$ for orbits $u_M$ and $u_{FS}$ and $y=0$ for $u_{R_h}$ and $u_{R_{pd}}$ --- see (\ref{eq:uM}), (\ref{eq:uFS}), (\ref{eq:uRH}) and (\ref{eq:uRPD}). We also reorder coordinates in the R\"ossler system (\ref{eq:rossler4D}) as $(y,x,z,w)$ so that in each case the section is defined by vanishing first coordinate.

With some abuse on notation, by the same letter $\mathcal P$ we will denote Poincar\'e map for each of these systems and by $u$ corresponding periodic orbit for this system. In each case $u$ is a fixed point of $\mathcal P^2$ (the orbit crosses section $\Pi$ twice in opposite directions). Let matrix $B$ be defined as in (\ref{eq:matricesVQB}).

Our aim is to compute $z = B^{-1}\left(\mathcal P^2(X)-u\right)$, where $X=u+B\intv{r}$, $\intv{r}=\frac{1}{2}s\cdot (0,[-1,1],\ldots,[-1,1])$ and $s>0$ is the diameter of the set $X$ in the coordinate system $B$. Observe, that set $X$ does not depend on the choice of vector $v_1$ used to define $V$ and $B$.

We compare three choices of matrix $A=B^{-1}$ and vector $y$ in Algorithm~\ref{alg:computePoincareMap}:

\begin{itemize}
    \item \textbf{cartesian} -- in this approach we compute Poincar\'e map in the coordinate system of the ODE by running Algorithm~\ref{alg:computePoincareMap} with arguments $A=\mathrm{Id}$, $y=0$ and then we change coordinates evaluating $z = B^{-1}\left(\mathcal P^2(X)-u\right)$,
    \item \textbf{diag+normal} --- we run Algorithm~\ref{alg:computePoincareMap} with its arguments $A=B^{-1}$, $y=u$, where $v_1$ orthogonal to the section,
    \item \textbf{diag+flowdir} --- we run Algorithm~\ref{alg:computePoincareMap} with its arguments $A=B^{-1}$, $y=u$, where $v_1=\frac{f(u)}{\|f(u)\|}$.
\end{itemize}
In all computation we have used a rigorous integrator from the CAPD library \cite{CAPDREVIEW,CAPD}. Since our goal is to check the properties of the algorithm that computes intersection of set of trajectories with section, we would like to reduce as much as possible the influence of rigorous integration along periodic trajectory. Therefore, we used the Hermite-Obreshkov method \cite{NedialkovJackson1998} of rather high order $50$ and high precision interval arithmetics with $200$ bits of mantissa \cite{MPFR}. In each case the initial condition was an affine set $X=u+B\intv{r}$, as described above, and represented as \texttt{Tripleton} initialized with  $X=(u, B, \intv{r}, \mathrm{Id}, \intv{0}, \mathrm{Id}, \intv{0})$.

\begin{figure}[htbp]
    \centerline{\includegraphics[height=.9\textheight]{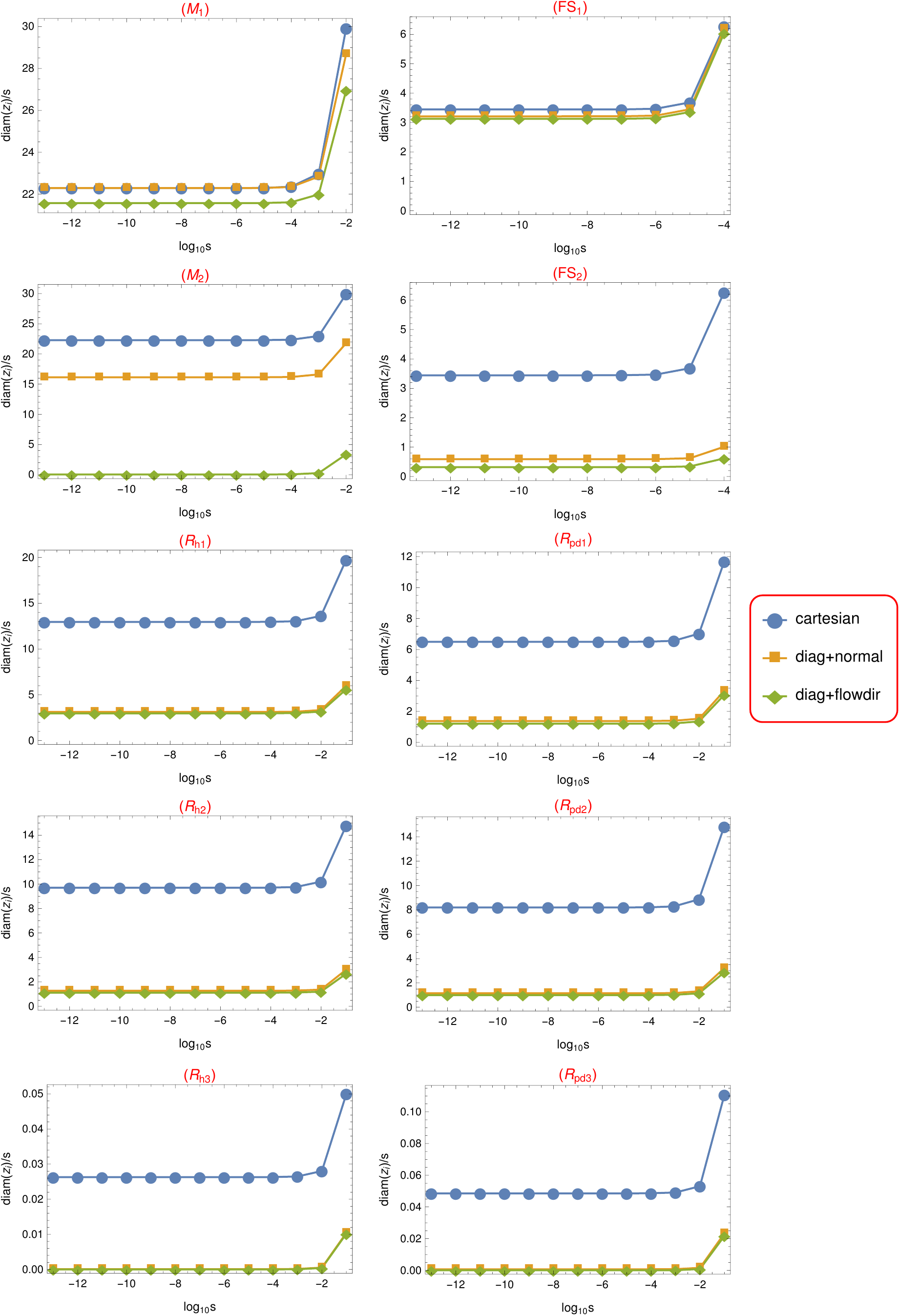}}
    \caption{Comparison of enclosures of $z = B^{-1}\left(\mathcal P^2(u+B\cdot \r)-u\right)$ computed by different methods. Horizontal axis corresponds to logarithms of size $s$ of the initial set $\r$. On vertical axis we plot the ratio $diam(z_i)/s$, which for small $s$ should be close to the corresponding eigenvalue of $D\mathcal P^2(u)$. The labels $M$, $FS$ and $R_h$ and $R_{pd}$ stand for the orbits $u_{M}$, $u_{FS}$, $u_{R_h}$ and $u_{R_{pd}}$, respectively. Additional index indicates coefficient of the result.\label{fig:fixedSection}}
\end{figure}

Results of experiments are partially listed in Tables~\ref{tab:michelson}, \ref{tab:fs}, \ref{tab:rossler} and \ref{tab:rosslerpd} (columns diag+normal and diag+flowdir) and visualised in Figure~\ref{fig:fixedSection}. These data strongly confirm that \textbf{diag+flowdir} strategy (that is Theorem~\ref{thm:sliding}) gives the best results. We also observe, that for small and moderate diameters of $X$, diameter of each component of the result is comparable to what is expected from linearization, that is scaling by corresponding eigenvalue --- see column diag+flowdir in Tables~\ref{tab:michelson}, \ref{tab:fs}, \ref{tab:rossler} and \ref{tab:rosslerpd}. This comment is not valid for the smallest in absolute value eigenvalue for the R\"ossler system, which is very close to zero and thus rounding errors and truncation errors of underlying ODE integrator are of comparable size.

\subsection{Poincar\'e section near periodic trajectories}\label{sec:variableSection}
In many applications we are free to choose Poincar\'e sections. Despite Theorem~\ref{thm:sliding}, which aims on reducing sliding effect on the section, we can try to minimize crossing time $\intv{\Delta t}$ by a proper choice of the section. This should additionally reduce errors coming from components $\intv{y}$ and $\intv{\Delta y}$ in Algorithm~\ref{alg:computePoincareMap}.

\begin{theorem}\label{thm:crossingTimeMin}
Let $f$ be a vector field and denote by $\varphi$ induced local flow on $\mathbb R^n$. Assume $\varphi(T,x_0)=x_0$ for some minimal $T>0$ and assume that $\lambda=1$ is an eigenvalue of $M:=D_x\varphi(T,x_0)$ of multiplicity one.

Let $\Pi=\{x\in\mathbb R^n\ :\ \alpha(x)=0\}$ be a local Poincar\'e section, where $\alpha:\mathbb R^n\to\mathbb R$ is smooth, $0$ is a regular value of $\alpha$ and $x_0\in \Pi$. Then
$$
\left(\ker Dt_{\Pi} = T_{x_0}\Pi\right) \Longleftrightarrow \left(D \alpha(x_0) \text{ is a left eigenvector of $M$ for $\lambda=1$}\right).
$$
In such case, $t_{\Pi}(x) = t_{\Pi}(x_0) + O(\|x-x_0\|^2)$ for $x\in\Pi$.
\end{theorem}
\begin{proof}
First we let us gather several easy facts which will be used in the proof. Observe that the fact that $\alpha$ is a local section we have
\begin{eqnarray}
  \scprod{ D\alpha(x_0) }{ f(x_0) } \neq 0,  \label{eq:trans-alpha} \\
  \ker D\alpha(x_0) = T_{x_0}\Pi.   \label{eq:kerDalpha}
\end{eqnarray}

    Differentiating the identity $\alpha(\varphi(t_\Pi(x),x))\equiv 0$ at $x=x_0$ we obtain
\begin{equation}\label{eq:gradId}
    \scprod{ D \alpha(x_0)}{ f(x_0)} D t_\Pi(x_0) + D\alpha(x_0)  M\equiv 0.
\end{equation}

It is well known that
\begin{equation}
  Mf(x_0) = f(x_0),  \label{eq:fx0-eigenvect}
\end{equation}
i.e.
 $f(x_0)$ in an eigenvector for $M$ for $\lambda=1$. From this and (\ref{eq:trans-alpha},\ref{eq:gradId}) we obtain
\begin{equation}
  \scprod{ D t_\Pi(x_0)}{ f(x_0) } = -\scprod{ D \alpha(x_0) }{ f(x_0) }^{-1}  \scprod{ D\alpha(x_0) }{ f(x_0)}\neq 0.  \label{eq:Dtf}
\end{equation}

Now we are ready for the proof of our equivalence.

For the proof in $\Leftarrow$ direction,  assume that $D\alpha(x_0)$ is a left eigenvector of $M$ for the eigenvalue $\lambda=1$. From (\ref{eq:gradId}) we obtain
\begin{eqnarray*}
  0=   \scprod{ D \alpha(x_0)}{ f(x_0)} D t_\Pi(x_0) +
    D\alpha(x_0) M =
    \scprod{ D \alpha(x_0) }{ f(x_0)} Dt_\Pi(x_0) + D\alpha(x_0).
\end{eqnarray*}
This implies that  $1$-forms $D \alpha(x_0)$ and $Dt_{\Pi}$ are proportional and from (\ref{eq:trans-alpha}) both are non-zero. Hence from
(\ref{eq:kerDalpha}) we obtain
\begin{equation*}
\ker D t_\Pi(x_0)=\ker D \alpha(x_0) =  T_{x_0}\Pi.
\end{equation*}

For the proof in $\Rightarrow$ direction assume  that $\ker D t_\Pi(x_0) = T_{x_0} \Pi$. From this and (\ref{eq:trans-alpha},\ref{eq:kerDalpha},\ref{eq:Dtf}) it follows that  the forms
$D \alpha(x_0)$ and $D t_\Pi(x_0)$ are non-zero and proportional. Hence from (\ref{eq:gradId}) we see that $D\alpha (x_0)$ is a left eigenvector of $M$ for some eigenvalue $\lambda$. We will show that $\lambda=1$.

We have
\begin{eqnarray*}
 \lambda  D \alpha(x_0) =  D \alpha(x_0) M,
\end{eqnarray*}
We evaluate both sides of the above equality in $f(x_0)$ and from (\ref{eq:fx0-eigenvect}) we obtain
\begin{eqnarray*}
  \lambda \scprod{ D \alpha(x_0)}{ f(x_0)} = \scprod{  D \alpha(x_0) M}{ f(x_0) } = \scprod{  D \alpha(x_0)}{ M f(x_0) }= \scprod{  D \alpha(x_0)}{ f(x_0)}
\end{eqnarray*}
Now from (\ref{eq:trans-alpha}) we infer that $\lambda=1$. This finishes the proof.
\end{proof}

In what follows, by CTO we will mean a \emph{crossing-time optimal} section chosen according to Theorem~\ref{thm:crossingTimeMin}.

\noindent\textbf{Example:} In order to show Theorem~\ref{thm:crossingTimeMin} in action let us consider the van der Pol equation
\begin{equation}
x' = y,\quad y' = 0.2y(1-x^2)-x. \label{eq:vanderpol}
\end{equation}
It is well known that the system has an attracting periodic point. A very accurate approximate initial condition for such periodic point is $u_0:=(x_0,y_0)=(2.0004136789920905,0.0)$

In the experiment, we will compute bounds on $t_{\mathcal P}(\intv{u})$ and $\mathcal P(\intv{u})$, where $\intv{u}=u_0 + \intv{\Delta u}$, for various sizes of the set $\intv{\Delta u}$ and for two different choices of the Poincar\'e section.

\begin{figure}[htbp]
\centerline{
    \includegraphics[width=.48\textwidth]{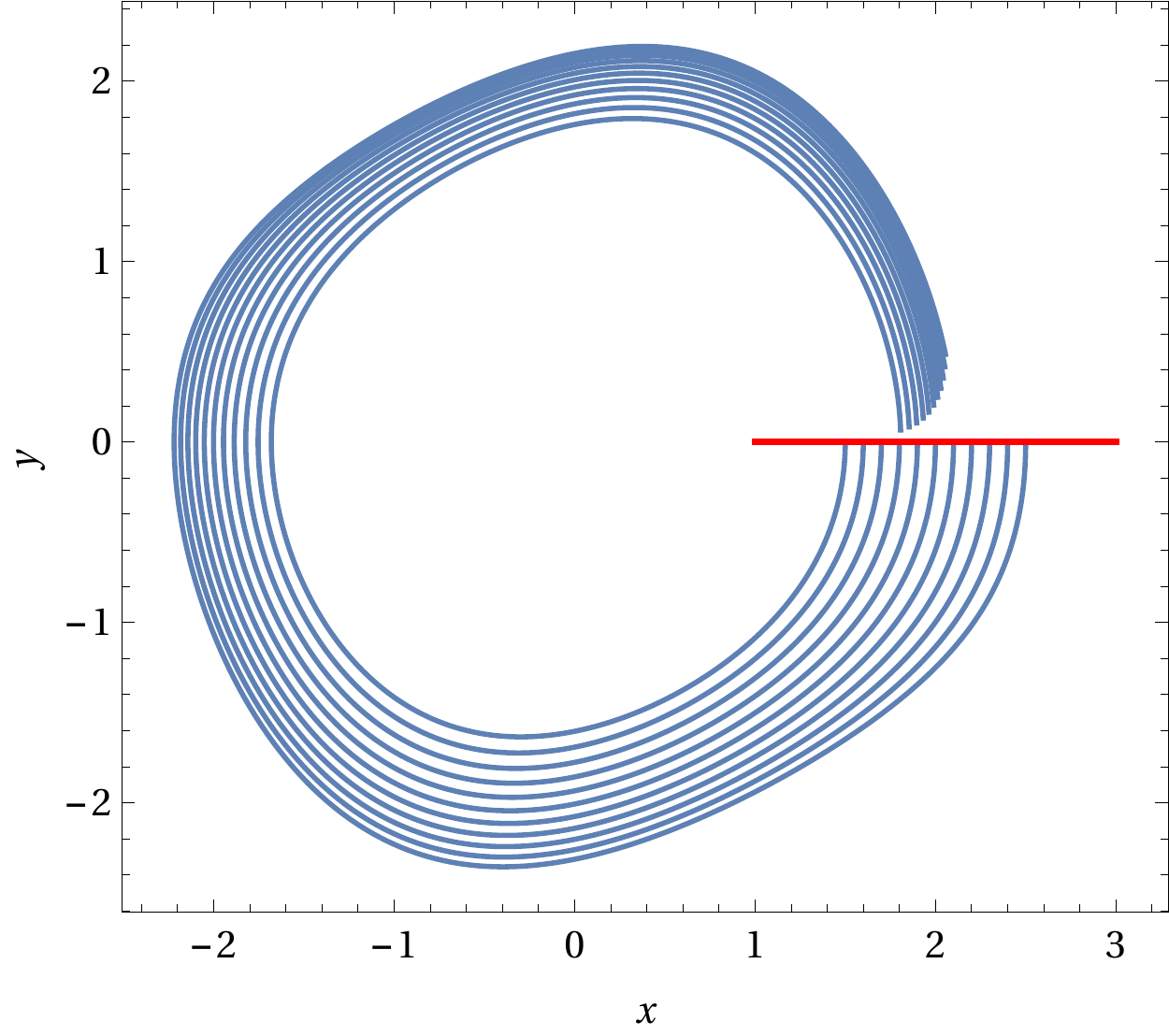}
    \includegraphics[width=.48\textwidth]{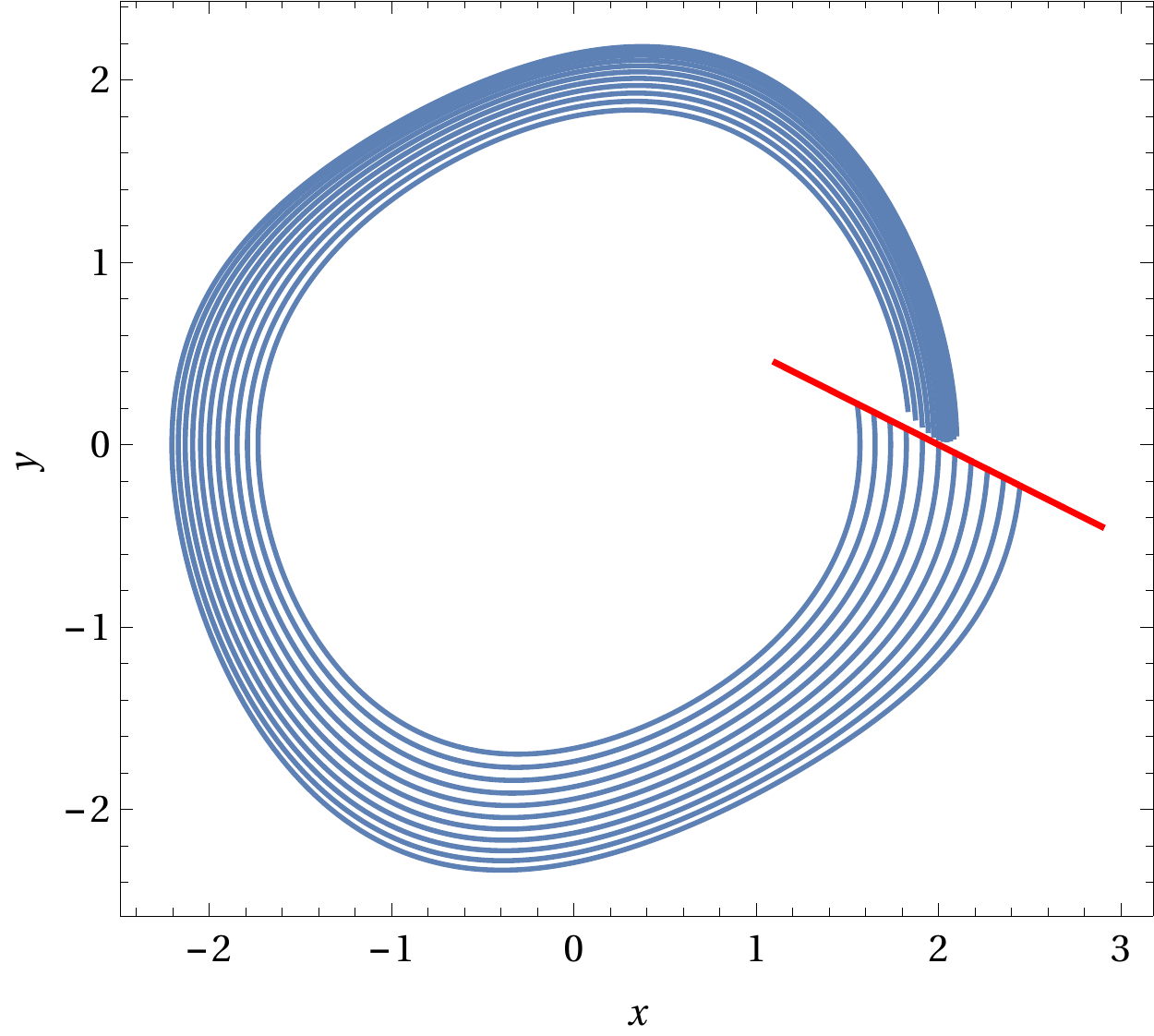}
}
\caption{Trajectory segments over the same time range of several points near an attracting periodic point for (\ref{eq:vanderpol}). Left: the Poincar\'e section is orthogonal to the flow at every point from the section. Right: Poincar\'e section chosen according to Theorem~\ref{thm:crossingTimeMin}, that is minimizing diameter of the crossing time. Here the return time has a local maximum at periodic point.\label{fig:crossingTime}}
\end{figure}

First, define $\Pi=\{(x,y)\in\mathbb R^2 : y=0\}$. This section is orthogonal to the vector field near the periodic point. In Figure~\ref{fig:crossingTime} (left panel) we have shown segments of trajectories of the same time-length for various initial conditions of the form $(x_0+\delta,y_0)$. It is a clear indication, that the return time as a function of $\delta$ is monotone near $\delta=0$. In Table~\ref{tab:monotoneCrossingTime} we give data from rigorous computation of Poincar\'e map obtained from Algorithm~\ref{alg:computePoincareMap} with coordinate system $A=\mathrm{Id}$ and $y=u_0$. We see that diameter of the result $\mathcal P(\intv{u})$ and of the crossing time grow almost linearly with the size of $\intv{\Delta u}= ([-\delta,\delta],0)$.

In Figure~\ref{fig:crossingTime} (right panel) we plot trajectory segments of the same time-length of several points at the CTO section near $u_0$. This is a numerical evidence, that the return time has a local maximum at periodic point. In Table~\ref{tab:flatCrossingTime} we give data from rigorous computation of Poincar\'e map obtained from Algorithm~\ref{alg:computePoincareMap} with a coordinate system $(x_1,x_2)$ centred at $u_0$ and given by normalized vector field at $u_0$, that is $(0,-1)$ and the vector spanning the section $(-0.894\ldots,0.449\ldots)$. As an initial condition for the Algorithm~\ref{alg:computePoincareMap} we choose sets of the form $\intv{u}=u_0+(0,[-\delta,\delta]x_2)$. We see that $x_2$ component grows linearly with $\delta$, as expected. Starting from $\delta=10^{-7}$ we see that diameter of return time is a flat function (quadratic growth) of $\delta$. For $\delta\leq 10^{-8}$ these values balance near $3.46\cdot 10^{-14}$ which is close to the machine epsilon in double precision.

\begin{table}[htbp]
    \begin{center}
        \small
\begin{tabular}{r|r|c}
$\delta = \frac{1}{2}\mathrm{diam}(\intv{u})$ & diameter of crossing time & $\pi_x\mathcal P(\intv{u})-x_0$\\ \hline
$10^{-9}$ & $3.6\cdot 10^{-10}$ & $[-2.83, 2.83]\cdot 10^{-10}$ \\ \hline
$10^{-8}$ & $3.6\cdot 10^{-9}$ & $[-2.83, 2.83]\cdot 10^{-9}$ \\ \hline
$10^{-7}$ & $3.6\cdot 10^{-8}$ & $[-2.83, 2.83]\cdot 10^{-8}$ \\ \hline
$10^{-6}$ & $3.6\cdot 10^{-7}$ & $[-2.83, 2.83]\cdot 10^{-7}$ \\ \hline
$10^{-5}$ & $3.6\cdot 10^{-6}$ & $[-2.83, 2.83]\cdot 10^{-6}$ \\ \hline
$10^{-4}$ & $3.61\cdot 10^{-5}$ & $[-2.83, 2.83]\cdot 10^{-5}$ \\ \hline
$10^{-3}$ & $3.64\cdot 10^{-4}$ & $[-2.84, 2.84]\cdot 10^{-4}$ \\ \hline
$10^{-2}$ & $3.97\cdot 10^{-3}$ & $[-2.93, 2.93]\cdot 10^{-3}$ \\ \hline
$10^{-1}$ & $1.18\cdot 10^{-1}$ & $[-6.5, 6.12]\cdot 10^{-2}$ \\ \hline

\end{tabular}
\end{center}

\caption{Computed diameter of crossing time and an enclosure of Poincar\'e map for the section $\Pi=\{(x,y) \ :\ y=0\}$ for (\ref{eq:vanderpol}) expressed in the cartesian coordinates.\label{tab:monotoneCrossingTime}}
\end{table}

\begin{table}[htbp]
    \begin{center}
        \small
        \begin{tabular}{r|r|c|c}
            $\delta = \frac{1}{2}\mathrm{diam}(\intv{u})$ & diameter of crossing time & $\pi_{x_2}\mathcal P(\intv{u})$\\ \hline
$10^{-9}$ & $3.46\cdot 10^{-14}$ & $[-2.83, 2.83]\cdot 10^{-10}$ \\ \hline
$10^{-8}$ & $3.46\cdot 10^{-14}$ & $[-2.83, 2.83]\cdot 10^{-9}$ \\ \hline
$10^{-7}$ & $6.39\cdot 10^{-14}$ & $[-2.83, 2.83]\cdot 10^{-8}$ \\ \hline
$10^{-6}$ & $2.99\cdot 10^{-12}$ & $[-2.83, 2.83]\cdot 10^{-7}$ \\ \hline
$10^{-5}$ & $2.96\cdot 10^{-10}$ & $[-2.83, 2.83]\cdot 10^{-6}$ \\ \hline
$10^{-4}$ & $2.96\cdot 10^{-8}$ & $[-2.83, 2.83]\cdot 10^{-5}$ \\ \hline
$10^{-3}$ & $2.97\cdot 10^{-6}$ & $[-2.83, 2.83]\cdot 10^{-4}$ \\ \hline
$10^{-2}$ & $3.11\cdot 10^{-4}$ & $[-2.89, 2.89]\cdot 10^{-3}$ \\ \hline
$10^{-1}$ & $6.26\cdot 10^{-2}$ & $[-4.66, 4.78]\cdot 10^{-2}$ \\ \hline

        \end{tabular}
    \end{center}
    \caption{Computed diameter of crossing time and an enclosure of Poincar\'e map for (\ref{eq:vanderpol}) for a CTO section expressed in $(x_1,x_2)$ coordinates.\label{tab:flatCrossingTime}}
\end{table}

\begin{remark}
Theorem~\ref{thm:crossingTimeMin} shows that using CTO section combined with Theorem~\ref{thm:sliding} near a periodic point we can reduce diameters of $\intv{y}$ and $\intv{\Delta y}$ to $diam(\intv{y}) \in O(diam(X_0)^3)$ and $diam(\intv{\Delta y)} \in O(diam(X_0)^4)$, because $\intv{\Delta t} \in O(diam(X_0)^2)$.
\end{remark}

\begin{figure}[htbp]
    \centerline{
        \includegraphics[width=.98\textwidth]{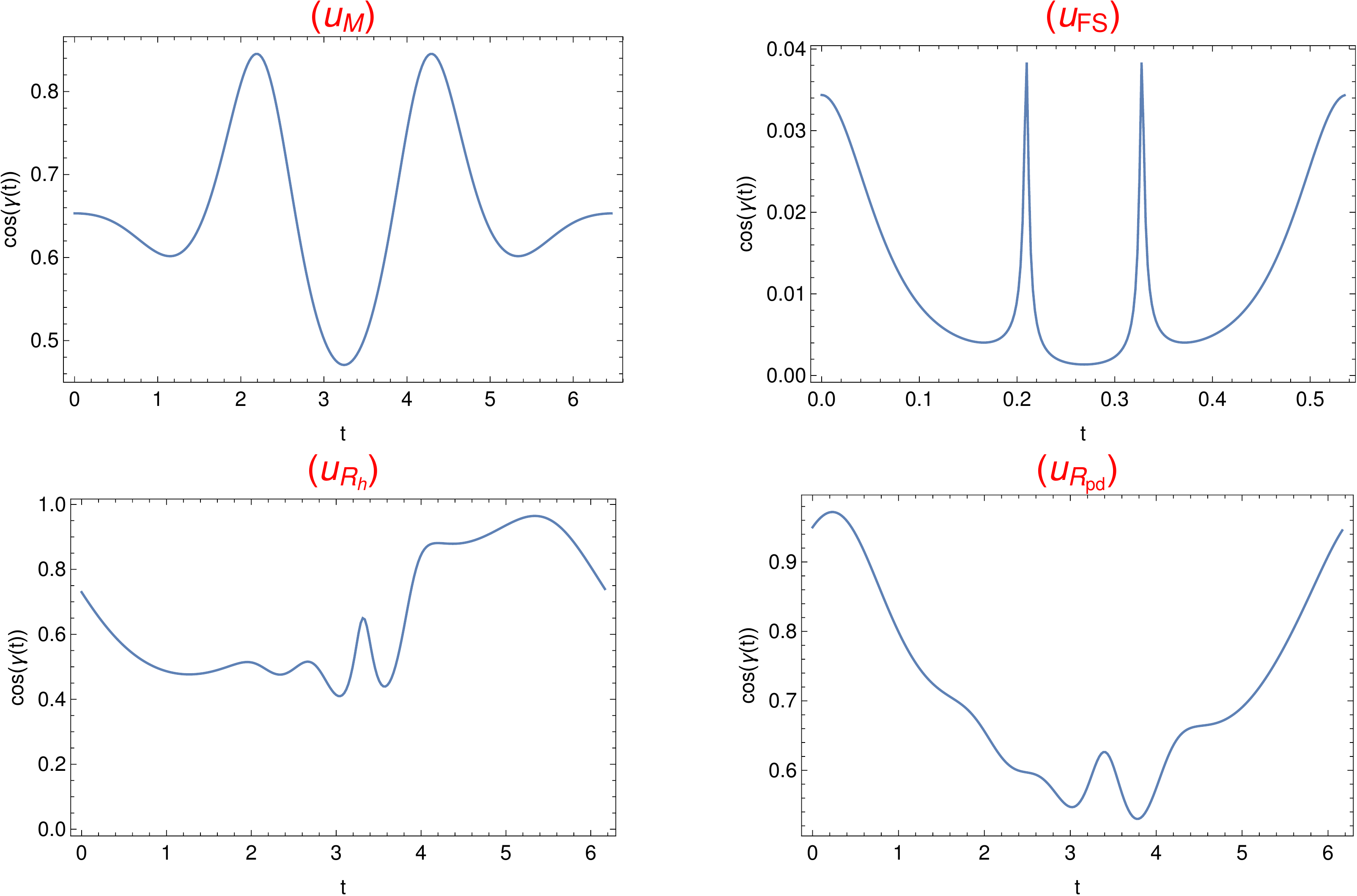}
    }
    \caption{Plot of $\cos(\gamma(t))$, where $\gamma(t)$ is the angle between vector field at $u(t)$ and normal vector to CTO section at $u(t)$ along periodic orbits $u_M$, $u_{FS}$, $u_{R_h}$ and $u_{R_{pd}}$, respectively. The plot $(u_{FS})$ shows that the vector field is nearly tangent to the CTO section along entire trajectory.\label{fig:angles}}
\end{figure}

\begin{remark}
    The CTO section may increase errors coming from sliding effect. In Figure~\ref{fig:angles} we plot the cosine of the angle between vector field at $u(t)$ and the normal vector to CTO section at $u(t)$ along orbits $u_M$, $u_{FS}$, $u_{R_h}$ and $u_{R_{pd}}$. We see that along entire trajectory of $u_{FS}$ the CTO section is almost tangent to the vector field.
\end{remark}

\subsection{Varying Poincar\'e section --- experiments.}
In this section we present results of the following experiment. We compare enclosures on Poincar\'e map for three choices of section for each periodic orbit $u_M$, $u_{FS}$, $u_{R_h}$ and $u_{R_{pd}}$, denoted below as $u$.
\begin{itemize}
    \item \textbf{orthogonal} --- the section is chosen as orthogonal to the vector field at $u$,
    \item \textbf{CTO} --- the section is chosen according to Theorem~\ref{thm:crossingTimeMin} at $u$,
    \item \textbf{max angle CTO} --- the section is chosen as CTO but at a point $u(t)$ on the trajectory in which the angle between vector field and normal to CTO section is the largest --- see Figure~\ref{fig:angles}.
\end{itemize}

Settings of the experiments are the same as in the case of fixed Poincar\'e section presented in Section~\ref{sec:fixedSections}. Our aim is to compute $z = B^{-1}\left(\mathcal P^2(X)-u\right)$, where $X=u+B\cdot \intv{r}$, $\intv{r}=\frac{1}{2}s\cdot (0,[-1,1],\ldots,[-1,1])$ and $s>0$ is the diameter of initial condition in the coordinate system $B$. In each case we run the Algorithm~\ref{alg:computePoincareMap} with the matrix $A=B^{-1}$ as described in \textbf{diag+flowdir} strategy (\ref{eq:matricesVQB}).

\begin{figure}[htbp]
    \centerline{\includegraphics[height=.93\textheight]{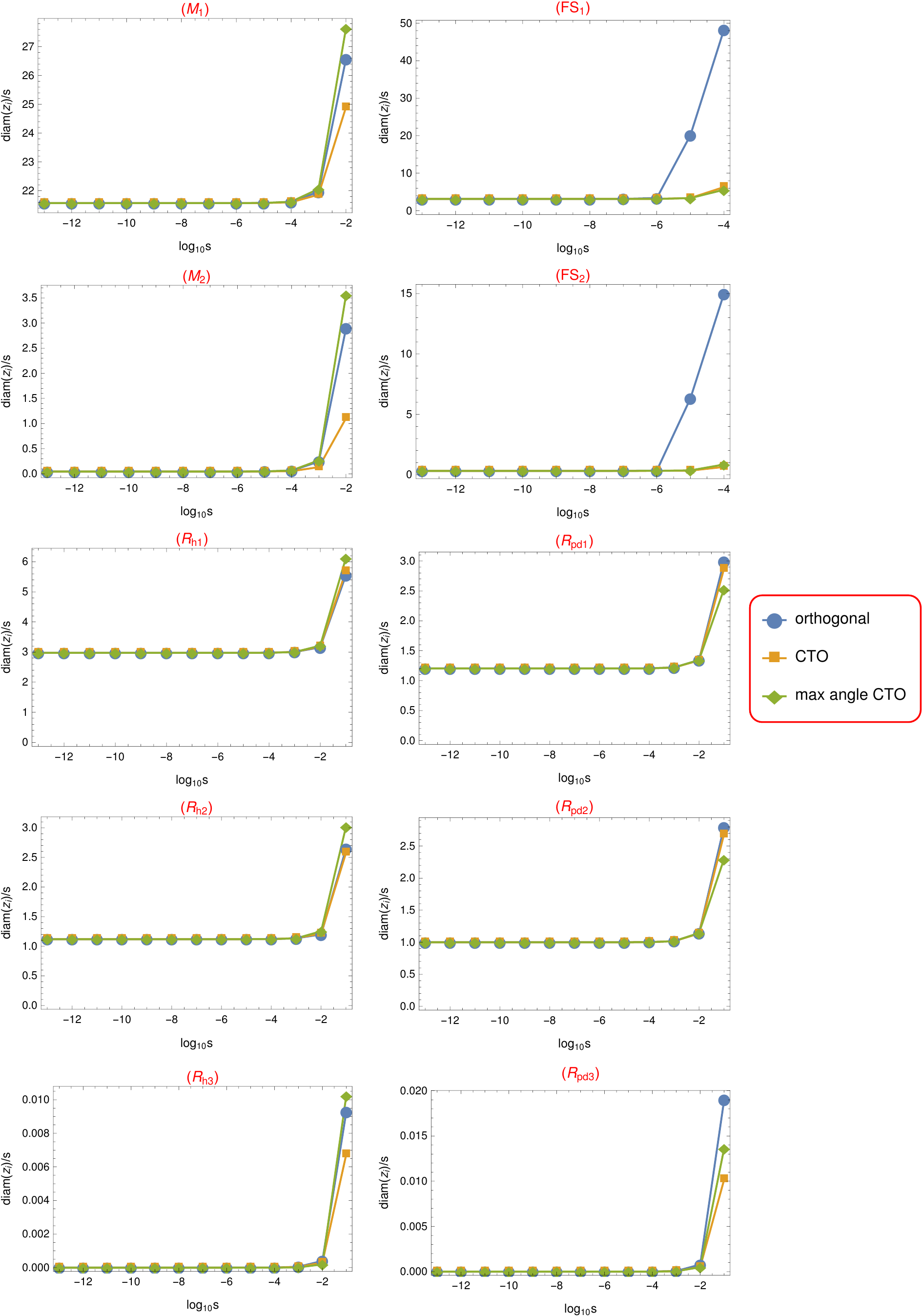}}
    \caption{Comparison of enclosures of $z = B^{-1}\left(\mathcal P^2(u+B\cdot \r)-u\right)$ computed by different methods. Horizontal axis corresponds to logarithms of size $s$ of the initial set $\r$. On vertical axis we plot the ratio $diam(z_i)/s$, which for small $s$ should be close to the corresponding eigenvalue of $D\mathcal P^2(u)$.\label{fig:varyingSection}}
\end{figure}

Results of experiment are partially listed in Tables~\ref{tab:michelson}, \ref{tab:fs}, \ref{tab:rossler} and \ref{tab:rosslerpd} (columns orthogonal, CTO and max angle CTO) and visualised in Figure~\ref{fig:varyingSection}. We would like to emphasize, that these data cannot be compared directly to each other and to results from previous experiment in Section~\ref{sec:fixedSections}. Here Poincar\'e maps have different domains and we actually integrate different sets (perhaps of different initial area/volume).

We observe that for small and moderate sizes of sets, results from experiments involving \textbf{diag+flowdir} strategy, that is fixed section, orthogonal, CTO and max angle CTO sections, are comparable. All of them are very accurate and there is no room for large improvements.

\begin{figure}[htbp]
    \centerline{\includegraphics[width=.75\textwidth]{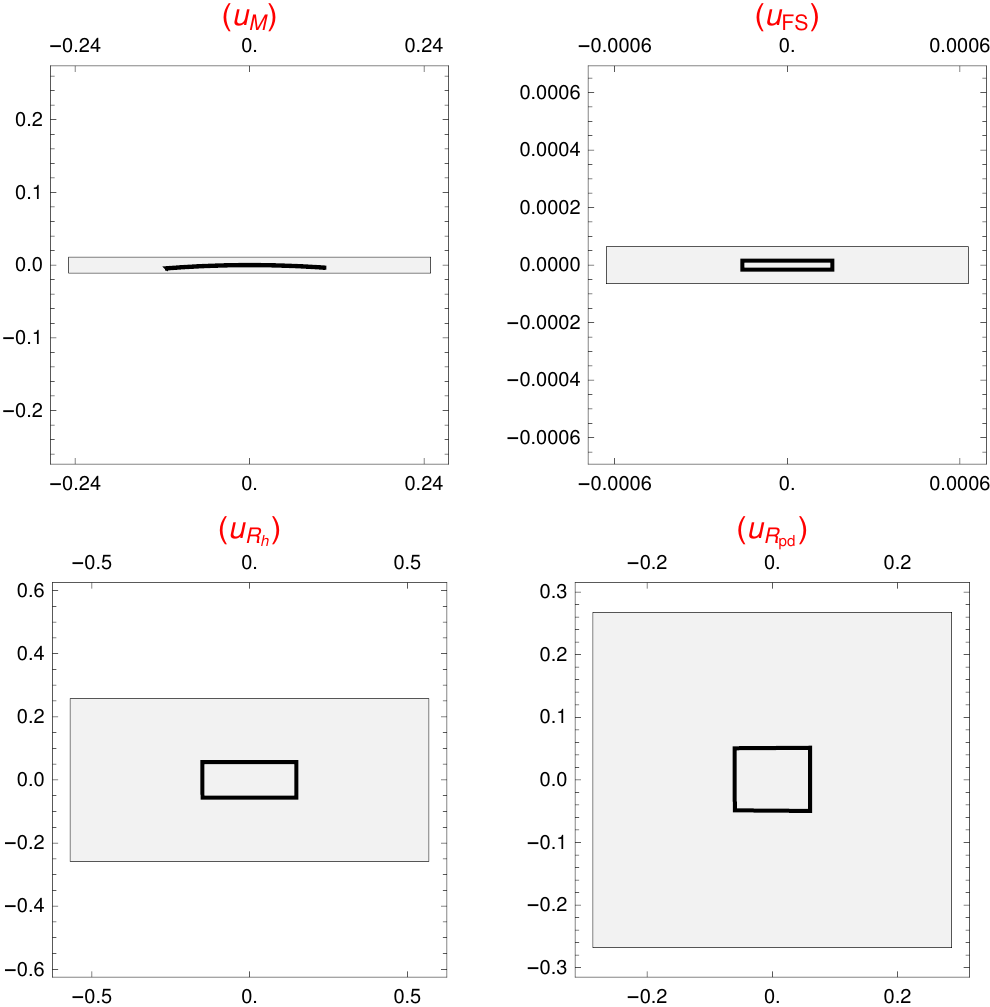}}
    \caption{In grey -- projection onto $(z_2,z_3)$ of an enclosure of $z = B^{-1}\left(\mathcal P^2(u+B\cdot \r)-u\right)$ onto two main coordinates computed for the largest size of $\intv{r}$ used in experiments --- see CTO column in Tables~\ref{tab:michelson}, Table~\ref{tab:fs}, Table~\ref{tab:rossler} and Table~\ref{tab:rosslerpd}. In black --- an approximate shape of such projection computed by a non-rigorous routine. In the case of 4D R\"ossler system we skip the last $z_4$ coordinate which corresponds to strong contraction.\label{fig:nonrig}}
\end{figure}

For large sets we see that the results are significantly overestimated --- see Figure~\ref{fig:nonrig}. The main source of such errors comes from integration of quite large set over the period of the orbit. This leads to overestimation of $X_0\supset \varphi(t_0,X)$ and in consequence to quite large enclosure on the main component of the result, that is $\intv{y_0}= \texttt{affineTransform}(X_0,A,y)$.

\begin{table}[htbp]
    \begin{center}
    \scriptsize
        \begin{tabular}{|c|c|c|c|c|c|}\hline
            $\log_{10}s$ &  \textbf{diag+normal} & \textbf{diag+flowdir} & \textbf{orthogonal} & \textbf{CTO} & \textbf{max angle CTO}\\ \hline
            \multicolumn{6}{|c|}{$\lambda_{M_1}\approx  -21.57189303583905$}\\ \hline
 -13 & 22.288770093234 & 21.571893035879 & 21.571893035877 & 21.571893035868 & 21.571893035886 \\
-12 & 22.288770093715 & 21.571893036244 & 21.571893036220 & 21.571893036127 & 21.571893036325 \\
-11 & 22.288770098524 & 21.571893039893 & 21.571893039664 & 21.571893038726 & 21.571893040654 \\
-10 & 22.288770146617 & 21.571893076394 & 21.571893074103 & 21.571893064721 & 21.571893083996 \\
-9 & 22.288770627515 & 21.571893441397 & 21.571893418486 & 21.571893324667 & 21.571893517420 \\
-8 & 22.288775436860 & 21.571897091434 & 21.571896862320 & 21.571895924126 & 21.571897851656 \\
-7 & 22.288823530421 & 21.571933591888 & 21.571931300742 & 21.571921918754 & 21.571941194109 \\
-6 & 22.289304477585 & 21.572298605637 & 21.572275692728 & 21.572181868412 & 21.572374627014 \\
-5 & 22.294115104218 & 21.575949664009 & 21.575720389491 & 21.574781702701 & 21.576709794306 \\
-4 & 22.342338639296 & 21.612552589485 & 21.610245257920 & 21.600813874537 & 21.620145477068 \\
-3 & 22.836610909239 & 21.988082379753 & 21.963506140123 & 21.864579188760 & 22.064023271359 \\
-2 & 28.684735372216 & 26.950444290822 & 26.572537443087 & 24.899391229219 & 27.640671971912 \\
\hline  \multicolumn{6}{|c|}{$ \lambda_{M_2}\approx -0.046356617768258279$}\\ \hline
 -13 & 16.134915640465 & 0.046356617787156 & 0.046356617784308 & 0.046356617774995 & 0.046356617788210 \\
-12 & 16.134915640887 & 0.046356617957266 & 0.046356617928783 & 0.046356617835662 & 0.046356617967803 \\
-11 & 16.134915645110 & 0.046356619658364 & 0.046356619373536 & 0.046356618442332 & 0.046356619763732 \\
-10 & 16.134915687348 & 0.046356636669347 & 0.046356633821066 & 0.046356624509024 & 0.046356637723029 \\
-9 & 16.134916109700 & 0.046356806779235 & 0.046356778296425 & 0.046356685175975 & 0.046356817316206 \\
-8 & 16.134920333587 & 0.046358507884187 & 0.046358223056044 & 0.046357291848356 & 0.046358613269082 \\
-7 & 16.134962572563 & 0.046375519540821 & 0.046372671255224 & 0.046363358859580 & 0.046376574908491 \\
-6 & 16.135384972877 & 0.046545696834452 & 0.046517213563835 & 0.046424057719320 & 0.046556402492512 \\
-5 & 16.139610031950 & 0.048253558313519 & 0.047968685063875 & 0.047033927584654 & 0.048375923353681 \\
-4 & 16.181967982389 & 0.065957180118245 & 0.063105333370683 & 0.053427383868315 & 0.067766038947423 \\
-3 & 16.616844425767 & 0.27530605365471 & 0.24486689683707 & 0.13335710379016 & 0.26868430398964 \\
-2 & 21.778676674834 & 3.3796717087139 & 2.8985584015404 & 1.1081551960197 & 3.5608490573157 \\\hline
        \end{tabular}
    \end{center}
    \caption{Computed ratio $diam(z_i)/s$ for the orbit $u_M$.
        Columns \textbf{diag+normal} and \textbf{diag+flowdir} contain results of the experiment with fixed Poincar\'e section as described Section~\ref{sec:fixedSections}. Columns orthogonal, CTO and max angle CTO contain results of the experiment with varying sections as described in Section~\ref{sec:variableSection}. Note that data in orthogonal, CTO and max angle CTO cannot be directly compared because these sections are different and thus we compute bound on Poincar\'e map on different sets.
        \label{tab:michelson}}
\end{table}

\begin{table}[htbp]
    \begin{center}
        \scriptsize
        \begin{tabular}{|c|c|c|c|c|c|}\hline
            $\log_{10}s$ &  \textbf{diag+normal} & \textbf{diag+flowdir} & \textbf{orthogonal} & \textbf{CTO} & \textbf{max angle CTO}\\ \hline
            \multicolumn{6}{|c|}{$\lambda_{FS_1}\approx  -3.1255162015308699$}\\ \hline
 -13 & 3.2116142444263 & 3.1255162038258 & 3.1255162269739 & 3.1255162039846 & 3.1255162077022 \\
-12 & 3.2116142655580 & 3.1255162244797 & 3.1255164559613 & 3.1255162260679 & 3.1255162249376 \\
-11 & 3.2116144768748 & 3.1255164310196 & 3.1255187458357 & 3.1255164469009 & 3.1255164321661 \\
-10 & 3.2116165900434 & 3.1255184964189 & 3.1255416446378 & 3.1255186552321 & 3.1255185078836 \\
-9 & 3.2116377217797 & 3.1255391504587 & 3.1257706385001 & 3.1255407385979 & 3.1255392650706 \\
-8 & 3.2118490440011 & 3.1257456956048 & 3.1280611613276 & 3.1257615776825 & 3.1257468381526 \\
-7 & 3.2139627521408 & 3.1278116218805 & 3.1510249035516 & 3.1279705113076 & 3.1278226902171 \\
-6 & 3.2351484960973 & 3.1485184344710 & 3.3866083502860 & 3.1501142089676 & 3.1485933519525 \\
-5 & 3.4519428524253 & 3.3604105774357 & 20.073719793613 & 3.3770717466865 & 3.3575308921411 \\
-4 & 6.1918831644386 & 6.0382345907118 & 48.192050322691 & 6.2929366399838 & 5.5793036988113 \\
\hline  \multicolumn{6}{|c|}{$ \lambda_{FS_2}\approx -0.31994714969328985$}\\ \hline
 -13 & 0.58904797067828 & 0.31994714992826 & 0.31994715230733 & 0.31994714994447 & 0.31994715050329 \\
-12 & 0.58904797368178 & 0.31994715204298 & 0.31994717583371 & 0.31994715220505 & 0.31994715370896 \\
-11 & 0.58904800371685 & 0.31994717319017 & 0.31994741109760 & 0.31994717481089 & 0.31994718904359 \\
-10 & 0.58904830406741 & 0.31994738466211 & 0.31994976374250 & 0.31994740086931 & 0.31994754319639 \\
-9 & 0.58905130758169 & 0.31994949938640 & 0.31997329079895 & 0.31994966145907 & 0.31995108473517 \\
-8 & 0.58908134341692 & 0.31997064711589 & 0.32020862211667 & 0.31997226791244 & 0.31998650120326 \\
-7 & 0.58938177100241 & 0.32018217307783 & 0.32256802070380 & 0.32019838802569 & 0.32034077392439 \\
-6 & 0.59239298018272 & 0.32230230643042 & 0.34678075911647 & 0.32246515567799 & 0.32389432028560 \\
-5 & 0.62320851091933 & 0.34399811694236 & 6.3011376294349 & 0.34569814280737 & 0.36052693835901 \\
-4 & 1.0129131821593 & 0.61828066926388 & 14.959060502378 & 0.64422908261053 & 0.84497758478011 \\\hline
        \end{tabular}
    \end{center}
    \caption{Computed ratio $diam(z_i)/s$ for the orbit $u_{FS}$.
        Columns \textbf{diag+normal} and \textbf{diag+flowdir} contain results of the experiment with fixed Poincar\'e section as described Section~\ref{sec:fixedSections}. Columns orthogonal, CTO and max angle CTO contain results of the experiment with varying sections as described in Section~\ref{sec:variableSection}. Note that data in orthogonal, CTO and max angle CTO cannot be directly compared because these sections are different and thus we compute bound on Poincar\'e map on different sets.\label{tab:fs}}
\end{table}

\begin{table}[htbp]
    \begin{center}
        \scriptsize
        \begin{tabular}{|c|c|c|c|c|c|}\hline
            $\log_{10}s$ &  \textbf{diag+normal} & \textbf{diag+flowdir} & \textbf{orthogonal} & \textbf{CTO} & \textbf{max angle CTO}\\ \hline
            \multicolumn{6}{|c|}{$\lambda_{R_{h1}}\approx  -2.9753618617896986$}\\ \hline
 -13 & 3.1269112296727 & 2.9753618617987 & 2.9753618617991 & 2.9753618617973 & 2.9753618617953 \\
-12 & 3.1269112296767 & 2.9753618618023 & 2.9753618618019 & 2.9753618618042 & 2.9753618618071 \\
-11 & 3.1269112299772 & 2.9753618620903 & 2.9753618620850 & 2.9753618620637 & 2.9753618620608 \\
-10 & 3.1269112327992 & 2.9753618647951 & 2.9753618647424 & 2.9753618645294 & 2.9753618645006 \\
-9 & 3.1269112610131 & 2.9753618918436 & 2.9753618913169 & 2.9753618891866 & 2.9753618888981 \\
-8 & 3.1269115432163 & 2.9753621623288 & 2.9753621570618 & 2.9753621357587 & 2.9753621328741 \\
-7 & 3.1269143652214 & 2.9753648671818 & 2.9753648145118 & 2.9753646014813 & 2.9753645726348 \\
-6 & 3.1269425854365 & 2.9753919158704 & 2.9753913891678 & 2.9753892588373 & 2.9753889703684 \\
-5 & 3.1272248040375 & 2.9756624185652 & 2.9756571512700 & 2.9756358454017 & 2.9756329603539 \\
-4 & 3.1300486358109 & 2.9783690272406 & 2.9783163272842 & 2.9781030120969 & 2.9780741257411 \\
-3 & 3.1551299308307 & 3.0024276270455 & 3.0019370205452 & 3.0053562752534 & 3.0006911559382 \\
-2 & 3.3312447246999 & 3.1612927762392 & 3.1579216074761 & 3.1824418580936 & 3.2139664668233 \\
-1 & 5.9748947614907 & 5.5919232064410 & 5.5499757238534 & 5.6828287198362 & 6.1202554505770 \\
         \hline  \multicolumn{6}{|c|}{$ \lambda_{R_{h2}}\approx 1.1193329361699592$}\\ \hline
 -13 & 1.2750081236836 & 1.1193329361772 & 1.1193329361780 & 1.1193329361755 & 1.1193329361756 \\
-12 & 1.2750081236814 & 1.1193329361751 & 1.1193329361749 & 1.1193329361760 & 1.1193329361771 \\
-11 & 1.2750081238407 & 1.1193329363308 & 1.1193329363269 & 1.1193329362917 & 1.1193329363020 \\
-10 & 1.2750081253204 & 1.1193329377781 & 1.1193329377391 & 1.1193329373876 & 1.1193329374899 \\
-9 & 1.2750081401109 & 1.1193329522510 & 1.1193329518609 & 1.1193329483465 & 1.1193329493692 \\
-8 & 1.2750082880882 & 1.1193330969806 & 1.1193330930796 & 1.1193330579353 & 1.1193330681626 \\
-7 & 1.2750097678299 & 1.1193345442775 & 1.1193345052677 & 1.1193341538239 & 1.1193342560970 \\
-6 & 1.2750245653503 & 1.1193490173467 & 1.1193486272464 & 1.1193451127757 & 1.1193461355106 \\
-5 & 1.2751725508514 & 1.1194937580732 & 1.1194898568592 & 1.1194547088893 & 1.1194649366501 \\
-4 & 1.2766534361840 & 1.1209421693982 & 1.1209031360826 & 1.1205513299068 & 1.1206536486960 \\
-3 & 1.2904231607513 & 1.1343998401486 & 1.1340331850761 & 1.1369430453153 & 1.1333056956081 \\
-2 & 1.3698919414549 & 1.2014426463159 & 1.1996643242106 & 1.2152067812507 & 1.2538154992768 \\
-1 & 2.9997723907998 & 2.6705679634689 & 2.6493196773279 & 2.5783568485068 & 3.0211083048782 \\\hline
    \end{tabular}
    \end{center}
    \caption{Computed ratio $diam(z_i)/s$ for the orbit $u_{R_h}$.
        Columns \textbf{diag+normal} and \textbf{diag+flowdir} contain results of the experiment with fixed Poincar\'e section as described Section~\ref{sec:fixedSections}. Columns orthogonal, CTO and max angle CTO contain results of the experiment with varying sections as described in Section~\ref{sec:variableSection}. Note that data in orthogonal, CTO and max angle CTO cannot be directly compared because these sections are different and thus we compute bound on Poincar\'e map on different sets.\label{tab:rossler}}
\end{table}

\begin{table}[htbp]
    \begin{center}
        \scriptsize
        \begin{tabular}{|c|c|c|c|c|c|}\hline 
            $\log_{10}s$ &  \textbf{diag+normal} & \textbf{diag+flowdir} & \textbf{orthogonal} & \textbf{CTO} & \textbf{max angle CTO}\\ \hline
            \multicolumn{6}{|c|}{$\lambda_{R_{pd1}}\approx  1.2039286263296685$}\\ \hline
 -13 & 1.3744854888330 & 1.2039286263323 & 1.2039286263322 & 1.2039286263322 & 1.2039286263316 \\
-12 & 1.3744854888497 & 1.2039286263484 & 1.2039286263480 & 1.2039286263474 & 1.2039286263453 \\
-11 & 1.3744854890396 & 1.2039286265322 & 1.2039286265287 & 1.2039286265228 & 1.2039286265000 \\
-10 & 1.3744854909230 & 1.2039286283550 & 1.2039286283202 & 1.2039286282607 & 1.2039286280334 \\
-9 & 1.3744855097563 & 1.2039286465824 & 1.2039286462353 & 1.2039286456404 & 1.2039286433675 \\
-8 & 1.3744856980897 & 1.2039288288571 & 1.2039288253857 & 1.2039288194367 & 1.2039287967079 \\
-7 & 1.3744875814245 & 1.2039306516056 & 1.2039306168909 & 1.2039305574015 & 1.2039303301127 \\
-6 & 1.3745064136666 & 1.2039488792229 & 1.2039485320727 & 1.2039479371714 & 1.2039456642602 \\
-5 & 1.3746947443819 & 1.2041311685904 & 1.2041276967474 & 1.2041217470111 & 1.2040990156302 \\
-4 & 1.3764909757232 & 1.2058681783762 & 1.2058331605620 & 1.2057750473739 & 1.2056335194386 \\
-3 & 1.3923461796892 & 1.2211702025606 & 1.2208174545532 & 1.2202191822506 & 1.2201485891365 \\
-2 & 1.5245370044750 & 1.3482398006188 & 1.3444079718526 & 1.3378502373169 & 1.3452772133533 \\
-1 & 3.3189752159528 & 3.0595461571246 & 2.9899841238967 & 2.8651359578762 & 2.5268610111985 \\
            \hline  \multicolumn{6}{|c|}{$ \lambda_{R_{pd2}}\approx -1$}\\ \hline
 -13 & 1.1486171722293 & 1.0000000000020 & 1.0000000000020 & 1.0000000000019 & 1.0000000000015 \\
-12 & 1.1486171722484 & 1.0000000000190 & 1.0000000000186 & 1.0000000000181 & 1.0000000000159 \\
-11 & 1.1486171724575 & 1.0000000002049 & 1.0000000002015 & 1.0000000001971 & 1.0000000001738 \\
-10 & 1.1486171745303 & 1.0000000020486 & 1.0000000020146 & 1.0000000019712 & 1.0000000017376 \\
-9 & 1.1486171952589 & 1.0000000204855 & 1.0000000201457 & 1.0000000197118 & 1.0000000173756 \\
-8 & 1.1486174025442 & 1.0000002048552 & 1.0000002014574 & 1.0000001971178 & 1.0000001737558 \\
-7 & 1.1486194753993 & 1.0000020485537 & 1.0000020145757 & 1.0000019711789 & 1.0000017375590 \\
-6 & 1.1486402030893 & 1.0000204856702 & 1.0000201458870 & 1.0000197119133 & 1.0000173756918 \\
-5 & 1.1488474901282 & 1.0002048700441 & 1.0002014718711 & 1.0001971314897 & 1.0001737670547 \\
-4 & 1.1508203440202 & 1.0019607036312 & 1.0019263987058 & 1.0018836409982 & 1.0017386848276 \\
-3 & 1.1681410842692 & 1.0174040964627 & 1.0170577983291 & 1.0165945694829 & 1.0165327959562 \\
-2 & 1.3108850142655 & 1.1452467277951 & 1.1414685208050 & 1.1359979815714 & 1.1438008444065 \\
-1 & 3.2096656391944 & 2.8618925282316 & 2.7920556484803 & 2.6772944104801 & 2.2936430842608 \\ \hline
        \end{tabular}
    \end{center}
    \caption{Computed ratio $diam(z_i)/s$ for the orbit $u_{R_{pd}}$.
        Columns \textbf{diag+normal} and \textbf{diag+flowdir} contain results of the experiment with fixed Poincar\'e section as described Section~\ref{sec:fixedSections}. Columns orthogonal, CTO and max angle CTO contain results of the experiment with varying sections as described in Section~\ref{sec:variableSection}. Note that data in orthogonal, CTO and max angle CTO cannot be directly compared because these sections are different and thus we compute bound on Poincar\'e map on different sets.\label{tab:rosslerpd}}
\end{table}

\section{Conclusions}
We proposed a new algorithm for efficient enclosing Poincar\'e maps. It takes advantage from the knowledge about representation of subsets of $\mathbb R^n$ by the underlying ODE solver. We have shown, that in general situation proper choice of coordinate system, that is \textbf{flowdir} strategy, allows to reduce sliding and wrapping effects when a set of trajectories crosses the section. Numerical experiments show, that for small and moderate sizes of sets, obtained enclosures are very close to what is expected from linear approximation of Poincar\'e map.

A strategy for choosing Poincar\'e sections at approximate periodic points, described in section \ref{sec:variableSection},  minimizes diameter of crossing time. This allows further reduction of overestimation of $\PM(X)$ for larger sets $X$. Although for large initial sets the main source of overestimation is a long time integration of ODE, our numerical findings in most of the cases indicate significant improvement given by CTO strategy when compared to the other ones.

The implementation of the algorithms presented in this paper are freely available as a part of the CAPD library \cite{CAPD,CAPDREVIEW}. A serious test for the quality of implementation was a computer--assisted verification of chaos in infinite dimensional system (Kuramoto--Sivashinsky PDE) \cite{WilczakZgliczynski2020} or verification of hyperchaos \cite{BarrioMartinezSerranoWilczak2015,WilczakSerranoBarrio2016} in the 4D R\"ossler system.

\bibliographystyle{elsarticle-num}
\addcontentsline{toc}{chapter}{Bibliography}
%uncomment next line to change bibliography name to references
\renewcommand{\bibname}{refs}
%\bibliographystyle{plain}  %use the plain bibliography style
%\section*{References}
\bibliography{refs}
\end{document}